\documentclass[12pt]{amsart}
\usepackage[centertags]{amsmath}

\usepackage{amsfonts}

\usepackage{amssymb}

\title[Indexed Hyperspaces]
{Acceptable Colorings of \\ Indexed Hyperspaces}
\author{James H. Schmerl}
\date{\today}

\def\into{\longrightarrow}

\def\harp{\hspace{-4pt} \upharpoonright \hspace{-4pt}}

\def\vR{\vec R}
\def\vS{\vec S}

\def\RR{{\mathbb R}}

\def\FF{{\mathbb F}}

 \DeclareMathOperator{\Th}{Th}

\DeclareMathOperator{\dom}{dom} 
\DeclareMathOperator{\supp}{supp} 
\DeclareMathOperator{\cov}{\tau} 
\DeclareMathOperator{\fcn}{fcn}

\begin{document}

\begin{abstract} Previous results about $n$-grids with acceptable colorings are extended here to 
$n$-indexed hyperspaces, which are structures ${\mathcal A} = (A;E_0,E_1, \ldots,E_{n-1})$,
where each $E_i$ is an equivalence relation on $A$.\end{abstract}

\maketitle

\bigskip

 If $1 \leq n < \omega$,  then, following \cite[Def.\@ 2.1]{sim97}, we say that   a 
structure ${\mathcal A} = (A;E_0,E_1, \ldots, E_{n-1})$ is an $n$-{\bf indexed hyperspace} if each $E_i$ is an equivalence relation on the set 
  $A$.   Given such an $n$-indexed hyperspace and $a \in A$, we let $[a]_i$ be the equivalence class of $E_i$ to which $a$ belongs.  A  {\bf coloring} of ${\mathcal A}$ is a function 
  $\chi : A \into n = \{0,1, \ldots, n-1\}$.
   The coloring $\chi$ is {\bf acceptable} 
 if whenever $a \in A$ and $i < n$, then 
the set $\{x \in [a]_i : \chi(x) = i\}$ is finite.  
 The Basic Question concerning   these notions is

\bigskip

\begin{center}

   {\em Which indexed hyperspaces have acceptable colorings}?
   
   \end{center}
\bigskip 

One of the incentives for considering this question is the still open instance of it concerning  
sprays. If $2 \leq m < \omega$ and $c \in \RR^m$ (where $\RR$ is the set of reals), then a {\bf spray centered at} $c$ is a set $S \subseteq \RR^m$ such that whenever $0 < r \in \RR$, then $\{x \in S : \|x-c\| = r\}$ is finite. The question 

\bigskip

\begin{center}

{\em How many sprays can cover $\RR^m$}?

\end{center}
\bigskip
was asked  in \cite[Question~2.4]{grid}. For  $m=2$, de la Vega \cite{dlv1}, answering an earlier question from \cite{sch03}, proved that 
$3$ sprays suffice to cover the plane $\RR^2$. 
In general, it follows from \cite{sik51} (or see Theorem~3.8) that it takes at least $m+1$ sprays to  cover $\RR^m$. On the other hand, 
as was observed in  \cite{grid},  it follows from \cite{sim97} (or see Theorem~3.2) that if $d < \omega$ and $2^{\aleph_0} \leq \aleph_d$, then $(d+1)(m-1)+1$ sprays do suffice to cover $\RR^m$. 

The questions about sprays can be reinterpreted into instances of the Basic Question. Given $c \in \RR^m$, let $E(c)$ be the equivalence relation on $\RR^m$ such that if $x,y \in \RR^m$, then $\langle x,y \rangle \in E(c)$ 
iff $\|x-c\| = \|y-c\|$. Then, for $c_0,c_1, \ldots,c_{n-1} \in \RR^m$, the $n$-indexed hyperspace 
$\big(\RR^m; E(c_0),E(c_1), \ldots, E(c_{n-1})\big)$ has an acceptable coloring iff there are sprays 
$S_0,S_1, \ldots, S_{n-1}$ centered at $c_0,c_1, \ldots, c_{n-1}$, respectively, such that 
$\RR^m = S_0 \cup S_1 \cup \cdots \cup S_{n-1}$. 

\bigskip

{\bf \S0.\@ Introduction.} An $n$-indexed hyperspace ${\mathcal A}$ is always understood to be so that ${\mathcal A} = (A; E_0,E_1, \ldots, E_{n-1})$. Some fundamental examples of $n$-indexed hyperspaces  are the $n$-cubes. 
 An $n$-indexed hyperspace  ${\mathcal A}$ is an 
 $n$-{\bf cube} if there are nonempty sets $A_0,A_1, \ldots, A_{n-1}$ such that 
 $A = A_0 \times A_1 \times \cdots \times A_{n-1}$ and whenever 
 $a,b \in A$ and $i < n$, then 
 $\langle a,b \rangle \in E_i$ iff $a_j = b_j$ for every $j < n$ such that $j \neq i$. 
 Thus, one can think of $[a]_i$ as ``the line through $a$ parallel to the $i$th coordinate axis.''
  We will call this ${\mathcal A}$ the $n$-cube {\bf for} $A$. For any set $X$, the $n$-cube {\bf over} $X$ is the $n$-cube for $X^n$.   The following classical theorem of Kuratowski answers the Basic Question for $n$-cubes over a set $X$. 

\bigskip


{\sc  Theorem 0.1}:  (Kuratowski \cite{kur51}) {\em Suppose that $1 \leq n < \omega$ and $X$ is set. Then, 
the $n$-cube over $X$ has an acceptable coloring iff $|X| < \aleph_{n-1} 
$.}


\bigskip

If the $n$-indexed hyperspace  ${\mathcal A}$ is such that 
$[a]_i \cap [a]_j$ is finite whenever $a \in A$ and $i < j < n$,
then (following \cite{grid}) we say that ${\mathcal A}$ is an $n$-{\bf grid}. Every  $n$-cube 
 is an $n$-grid.  
 
 The right-to-left half of Kuratowski's Theorem~0.1 is a consequence of the following  more general theorem, which, itself, is a consequence of the still more general \cite[Theorems~5.1 \& 5.2]{grid}.

 \bigskip
 
 
 {\sc Theorem 0.2}: {\em If $1 \leq n < \omega$, ${\mathcal A}$ is an $n$-grid and $|A| < \aleph_{n-1}$, then ${\mathcal A}$ 
 has an acceptable coloring.}
 
 
 \bigskip
 
If  ${\mathcal A}$ and ${\mathcal B}$ are $n$-indexed hyperspaces, then an {\bf  embedding} of 
${\mathcal B}$ into ${\mathcal A}$ is defined, as expected, to be  a one-to-one function $f : B \into A$ such that  whenever $x,y \in B$ and $i < n$, then
$$
[x]_i = [y]_i \Longleftrightarrow  [f(x)]_i = [f(y)]_i.
$$
If there is an embedding of ${\mathcal B}$ into ${\mathcal A}$, then we say that 
${\mathcal B}$ is  {\bf embeddable} into ${\mathcal A}$ or that 
${\mathcal A}$  {\bf embeds} ${\mathcal B}$. Obviously, if ${\mathcal A}$ embeds ${\mathcal B}$ and 
 ${\mathcal A}$ has an acceptable coloring, then ${\mathcal B}$ has an  acceptable coloring.

A consequence of the left-to-right half of Kuratowski's Theorem is that if $|X| \geq \aleph_{n-1}$ and 
${\mathcal A}$ is an $n$-indexed hyperspace that embeds the $n$-cube over $X$, then ${\mathcal A}$ does not have an acceptable coloring. De la Vega proved a  partial converse to this for $n$-grids.

\bigskip


{\sc Theorem 0.3}: (de la Vega \cite{dlv2}) {\em  If ${\mathcal A}$ is an $n$-grid that does not embed every finite 
$n$-cube, then ${\mathcal A}$ has an acceptable coloring.}


\bigskip

The converse of de la Vega's Theorem is not true in general (as Kuratowski's Theorem shows). There are even arbitrarily large $n$-grids that embed every finite $n$-cube and  have acceptable colorings. However,  Theorem~0.2 is the only obstacle to the converse of de la Vega's Theorem 
 when it is restricted to  semialgebraic grids (a definition of which is given in \S5).   Thus, the Basic Question is answered for semialgebraic grids by Theorem~0.2 and the following theorem \cite[Lemma~3.6 \& Coro.\@ 4.3]{grid}.

\bigskip


{\sc Theorem 0.4:} {\em Suppose that ${\mathcal A}$ is a semialgebraic $n$-grid and $|A| \geq \aleph_{n-1}$.  The following are equivalent$:$

$(1)$ ${\mathcal A}$ has an acceptable coloring.

$(2)$ ${\mathcal A}$ does not embed every finite 
$n$-cube.

$(3)$ ${\mathcal A}$ does not embed the $n$-cube over  $\RR$.} 


\bigskip

A consequence of Theorem~0.2 and the proof \cite{grid} of Theorem~0.4 is the following theorem concerning decidability.
See \S5 for more of an explanation and also for a generalization to indexed hyperspaces.

\bigskip

{\sc Theorem 0.5}: {\em The set of  ${\mathcal L}_{OF}$-formulas  that, for some $n < \omega$, define  a semialgebraic $n$-grid having an acceptable coloring is computable.}

\bigskip

The Basic Question for grids was studied in \cite{grid}. 
Our aim in this paper is to extend results about acceptable colorings of $n$-grids to $n$-indexed hyperspaces. We will do so for all the results of \cite{grid}. 

The outline of the rest of this paper is as follows. The easy answer to the Basic Question for countable 
indexed hyperspaces is given in \S1. A characterization of $n$-grids having acceptable colorings was given by de la Vega in \cite{dlv1} and \cite{dlv2}. A characterization for 
indexed hyperspaces in the spirit of de la Vega's is in \S2. An important  step in generalizing 
results about $n$-grids to $n$-indexed hyperspaces was already undertaken by  Simms \cite{sim97}.  His generalization of 
Theorems~0.2  is discussed in   \S3 but in a way that differs from what is in \cite{sim97}.  
That section also contains an improvement and simplification of his generalization of Theorem~0.1.

Theorem~0.3 will also be extended to indexed hyperspaces in Theorem~4.2. Even when Theorem~4.2 is restricted to grids (Corollary~4.3), this yields an improvement of Theorem~0.3. 
The extension of Theorem~0.3 to indexed hyperspaces is presented and proved in \S4.  
A strengthened version of Theorem~0.4 is proved in \cite{grid}, yielding Theorem~0.5 as a consequence. These results will be extended to semialgebraic indexed hyperspaces in \S5, yielding   the decidability of the set of formulas defining  semialgebraic indexed hyperspaces having acceptable colorings.
Thus, in principle, the question of how many sprays are needed to cover $\RR^n$  should be answerable. 

 \bigskip

{\bf \S1.\@ Countable indexed hyperspaces.} The main result of this short section, Corollary~1.3, 
 characterizes  those  countable $n$-indexed hyperspaces that have acceptable colorings. We start 
 with a simple lemma in which there is no countability condition.
 
 \bigskip
 
 
 {\sc Lemma 1.1}: {\em Suppose that ${\mathcal A}$ is an $n$-indexed hyperspace that  has an acceptable coloring. Then, for every $a \in A$,  $[a]_0 \cap [a]_1 \cap \cdots \cap [a]_{n-1}$  is finite.}
 

\bigskip

{\it Proof}. 
Suppose that $a \in A$ and $X = [a]_0 \cap [a]_1 \cap \cdots \cap [a]_{n-1}$ is infinite. Let $\chi : A \into n$ be a coloring. By the Pigeon Hole Principle, there is $i < n$ such that 
$X \cap \chi^{-1}(i)$ is infinite. Then, $\{x \in [a]_i : \chi(x) = i\}$ is infinite, so $\chi$ is not acceptable.
\qed

\bigskip

Next, we show that the converse of Lemma~1.1 holds when restricted to countable ${\mathcal A}$.

\bigskip


{\sc Lemma 1.2:} {\em Suppose that ${\mathcal A}$  is a countable $n$-indexed hyperspace.   If $[a]_0 \cap [a]_1 \cap \cdots \cap [a]_{n-1}$  is finite for all $a \in A$, then ${\mathcal A}$ has an acceptable coloring.}


\bigskip

{\it Proof}.    
Suppose that   $[a]_0 \cap [a]_1 \cap \cdots \cap [a]_{n-1}$  is finite whenever $a \in A$.
We can assume that $A$ is infinite, as otherwise every coloring is acceptable. 
Hence, let $a_0,a_1,a_2, \ldots$ be a nonrepeating enumeration of~$A$. 

To define $\chi : A \into n$, consider an arbitrary $a = a_k \in A$. For each $i < n$, let $m_i$ be the least $m < \omega$ such that  $a \in [a_m]_i$. Notice that each $m_i$ is well-defined and that $m_i \leq k$.
 Then let $\chi(a) = j$, where $m_j = \max\{m_i : i < n\}$. (Since there may be more than one possible such $j$, to be definitive, choose the least one.)  This defines $\chi$, which clearly  
 is a coloring of~${\mathcal A}$.
 
 We claim that $\chi$ is acceptable. For a contradiction, suppose that $j < n$, $a \in A$, $X \subseteq [a]_j$ 
 is infinite and $\chi$ is constantly $j$ on $X$. By the maximality in the definition of $\chi$, for each $i < n$ and $x \in X$, there is $r \leq j$ such that $x \in [a_r]_i$. By the Pigeon Hole Principle, we can assume that there are 
 $r_0,r_1, \ldots, r_{n-1} \leq j$ such that $X \subseteq [a_{r_i}]_i$ for each   $i < n$. But then, taking $a \in X$, we have that 
 $X \subseteq [a]_0  \cap [a]_1 \cap  \cdots \cap [a]_{n-1}$. Since $X$ is infinite, this contradicts our assumption, thereby proving that $\chi$ is acceptable.   \qed
 
 \bigskip
 
 
 {\sc Corollary 1.3}: {\em Suppose that ${\mathcal A}$  is a countable $n$-indexed hyperspace. 
  Then, ${\mathcal A}$ has an acceptable coloring iff $[a]_0 \cap [a]_1 \cap \cdots \cap [a]_{n-1}$  is finite for all $a \in A$.} \qed
  
  
  \bigskip


{\bf \S2.\@  Twisted  indexed hyperspaces.}  
Using elementary substructures of the set-theoretic universe, de la Vega \cite{dlv1} defined the notion of 
a twisted $3$-grid and  proved that  a $3$-grid is twisted iff it has an acceptable coloring. 
Later \cite{dlv2}, he extended the definition to all $n$-grids and proved that an $n$-grid is twisted iff it has an acceptable coloring. We define here a closely related notion that is applicable to all $n$-indexed hyperspaces. This definition uses 
only elementary substructures of the indexed hyperspace, but an approach closer to  de la Vega's  would work just as well.    Since the consequences, at least for $n$-grids, are the same, 
we have decided to appropriate 
 de la Vega's  term in Definition~2.1. The main result of this section is Theorem~2.2.
 
 We will need a minor  generalization of terminology. 
 If $I \subseteq \omega$ is finite, then 
${\mathcal A}$ is 
an $I$-{\bf indexed hyperspace} if ${\mathcal A} = (A; E_i)_{i \in I}$, where each $E_i$ is an 
equivalence relation on $A$. If ${\mathcal A}$ is such 
an $I$-indexed hyperspace and $J \subseteq I$, then we let ${\mathcal A} \harp J$ be the 
$J$-indexed hyperspace  $ (A; E_j)_{j \in J}$.

 Suppose, for the moment,  that  ${\mathcal A} = (A; \ldots)$ is any first-order structure. If $B \subseteq A$, then we let  ${\mathcal A}|B$ be the substructure of ${\mathcal A}$ with universe $B$ (if there is such a substructure). If $|A| = \kappa > \aleph_0$, then we define 
a {\bf filtration} for ${\mathcal A}$ to be 
a sequence $\langle A_\alpha : \alpha < \kappa \rangle$ of subsets of $A$ 
such that $|A_\alpha| < \kappa$ for each $\alpha < \kappa$  and $\langle {\mathcal A}|A_\alpha :  \alpha < \kappa \rangle$ is an increasing, continuous chain of  elementary substructures 
of ${\mathcal A}$ whose union is ${\mathcal A}$.   Every  ${\mathcal A}$ (for a countable language) of uncountable cardinality $\kappa$ has a filtration $\langle A_\alpha: \alpha < \kappa \rangle$ with the additional property that each $|A_\alpha| = |\alpha| + \aleph_0$.

\bigskip


{\sc Definition} 2.1: (by recursion) Suppose that ${\mathcal A}$ is an $n$-indexed hyperspace 
and $|A| = \kappa$. 
We say that  ${\mathcal A}$ is {\bf twisted} if  

\smallskip

$(0)$ $[a]_0 \cap [a]_1 \cap \cdots \cap [a]_{n-1}$ is finite whenever $a \in A$, 

\smallskip

\noindent and either

\smallskip

$(1)$  ${\mathcal A}$ is countable, 

\smallskip

\noindent or else

\smallskip

$(2)$ $|A| = \kappa > \aleph_0$ and there is a filtration $\langle A_\alpha : \alpha < \kappa \rangle$ for ${\mathcal A}$ such that ${\mathcal A}|A_0$ is twisted and whenever $\alpha < \kappa$,
$\varnothing \neq I \subseteq n$ and 
$$
B = \{x \in A_{\alpha+1} \backslash A_\alpha : \forall i < n[ i \in I \leftrightarrow [x]_i \cap A_\alpha =\varnothing ]\},
$$ 
then $({\mathcal A}|B)\harp I$ is twisted.


\bigskip

If  ${\mathcal A}$ is an uncountable $n$-indexed hyperspace, then we will refer to a filtration for  
${\mathcal A}$ as in (2) of Definition~2.1 as a {\bf twisted} filtration. 

\bigskip


{\sc Theorem 2.2:} {\em If ${\mathcal A}$ is an $n$-indexed hyperspace, then ${\mathcal A}$ is twisted iff it has an acceptable coloring.}


\bigskip

{\it Proof}. The theorem will be proved by induction on the cardinality of~$A$.  Corollary~1.3 proves the theorem in case ${\mathcal A}$ is countable. Now assume that $|A| = \kappa > \aleph_0$ and that the theorem is true for all smaller indexed hyperspaces. 
\smallskip

$(\Longrightarrow)$:  Suppose that ${\mathcal A}$ is twisted. Let $\langle A_\alpha : \alpha < \kappa \rangle$ be a twisted flirtation for ${\mathcal A}$, and let ${\mathcal A}_\alpha = {\mathcal A}|A_\alpha$ for $\alpha < \kappa$.  We will obtain, by transfinite recursion, 
a sequence $\langle \chi_\alpha : \alpha < \kappa \rangle$ such that whenever $\alpha < \beta < \kappa$, then:

\begin{itemize}

\item $\chi_\alpha$ is an acceptable coloring of ${\mathcal A}_\alpha$;

\item $\chi_\alpha \subseteq \chi_\beta$;

\item if  $x \in A_\beta \backslash A_\alpha$ and $\chi_\beta(x) = i$, then  $[x]_i \cap A_\alpha = \varnothing$.

\end{itemize}
We then will have that $\chi = \bigcup_{\alpha < \kappa}\chi_\alpha$ is an acceptable coloring of ${\mathcal A}$.

Since $|A_0| < \kappa$ and ${\mathcal A}_0$ is twisted,  then,  by the inductive hypothesis, ${\mathcal A}_0$  has an acceptable coloring  $\chi_0$.  

If $\alpha$ is a limit ordinal, then let $\chi_\alpha = \bigcup_{\gamma< \alpha}\chi_\gamma$.

We now come to the case of successor ordinals. Suppose that we have $\chi_\gamma$ for $\gamma \leq \alpha$. For each $I \subseteq n$, let $B_I$ be defined as $B$ is in Definition~2.1(2). 

We will show that $B_{\varnothing} = \varnothing$. To the contrary,  suppose that $x \in B_{\varnothing}$. 
Then, $[x]_i \cap A_\alpha \neq \varnothing$ for each $i < n$. Let $y_i \in [x]_i \cap A_\alpha$ for each $i < n$. Then  $x \in [y_0]_0 \cap [y_1]_1 \cap \cdots \cap [y_{n-1}]_{n-1}$  so that $[y_0]_0 \cap [y_1]_1 \cap \cdots \cap [y_{n-1}]_{n-1}$ has a nonempty intersection with $A_{\alpha +1} \backslash A_\alpha$. Then, by elementarity, 
$[y_0]_0 \cap [y_1]_1 \cap \cdots \cap [y_{n-1}]_{n-1}$ has an infinite intersection $D$ with $A_\alpha$. For any $a \in D$, we have that $[a]_0 \cap [a]_1 \cap \cdots \cap [a]_{n-1} \supseteq D$, 
contradicting~(0). 
 
Therefore, for each $x \in A_{\alpha+1} \backslash A_{\alpha}$, there is a unique nonempty 
$I \subseteq n$ such that $x \in B_I$. 

If $\varnothing \neq I \subseteq n$, then  ${\mathcal B}_I = ({\mathcal A}|B_I)\harp I$ is twisted, so, by the inductive 
 hypothesis, we can let $\varphi_I : B_I \into I$ be an acceptable coloring of ${\mathcal B}_I$. 
 Then let $\chi_{\alpha+1} = \chi_\alpha \cup \bigcup_{\varnothing \neq I\subseteq n}\varphi_I$.

\smallskip 

$(\Longleftarrow)$: Lemma~1.1 shows that $(0)$ holds whenever ${\mathcal A}$ has an acceptable coloring. This takes care of the case  of countable indexed hyperspaces. 
For uncountable ones, we will prove the following by induction on $\kappa$:

\begin{quote}

Suppose that  ${\mathcal A}$ is an indexed hyperspace, $\chi$ is an acceptable coloring of ${\mathcal A}$ and $|A| = \kappa > \aleph_0$. Then every filtration for $({\mathcal A}, \chi)$ is a twisted filtration 
for ${\mathcal A}$. 

\end{quote}

 Let 
$\langle A_\alpha : \alpha < \kappa \rangle$ be a filtration for the expanded structure 
$({\mathcal A}, \chi)$, and let ${\mathcal A}_\alpha = {\mathcal A}|A_\alpha$ for $\alpha < \kappa$.

First, we show that ${\mathcal A}_0$ is twisted. If $A_0$ is countable, then ${\mathcal A}_0$ is twisted 
since, by Lemma~1.1, ${\mathcal A}$ and, consequently, ${\mathcal A}_0$ satisfy $(0)$. 
If $A_0$ is uncountable, then, since ${\mathcal A}$ and, consequently, ${\mathcal A}_0$ have acceptable colorings, then, by the inductive hypothesis, ${\mathcal A}_0$ is twisted.

 Next, consider $\alpha < \kappa$ and nonempty $I \subseteq n$. Define  $B$  as in Definition~2.1(2), 
 and let ${\mathcal B} = ({\mathcal A}|B)\harp I$. We want to show that ${\mathcal B}$ is twisted. 

To prove that ${\mathcal B}$ is twisted, it suffices to prove that it has an acceptable coloring. 
We will do so by showing that, in fact, $\chi \harp B$ is an acceptable coloring of ${\mathcal B}$. If $\chi \harp B$ is a coloring, then clearly it is acceptable, so we need only show that $\chi \harp B$ is a coloring. Let $x \in B$ and suppose, for a contradiction, that $\chi(x) = i \not\in I$. That implies that $[x]_i \cap  A_\alpha \neq \varnothing$. Since $[x]_i \cap \chi^{-1}(i)$ 
is finite, it follows by elementarity that $[x]_i \cap \chi^{-1}(i) \subseteq A_\alpha$, so $x \in A_\alpha$, which is a contradiction. \qed 

\bigskip

If $X$ is any set, then ${\mathcal P}(X)$ is the set of subsets of $X$, and if $n< \omega$, then
$[X]^n$ is the set of $n$-element subsets of $X$.

Suppose that $n < \omega$ and ${\mathcal I} \subseteq {\mathcal P}(n)$.  We say that 
              ${\mathcal A}$ is an $(n, {\mathcal I})$-{\bf grid} if it is an $n$-indexed hyperspace such that  whenever $I \in {\mathcal I}$ and $a \in A$,
              then $\bigcap_{i \in I}[a]_i$ is finite. 
              
                     We present some examples of  $(n,{\mathcal I})$-grids. Suppose that ${\mathcal A}$ is an $n$-indexed hyperspace. Vacuously, ${\mathcal A}$ is an $(n,\varnothing)$-grid, and, conventionally, ${\mathcal A}$ is finite  iff it is an $(n,\{\varnothing\})$-grid. By Lemma~1.1, if ${\mathcal A}$ has an acceptable coloring, then ${\mathcal A}$ is an $(n,\{n\})$-grid.  Lastly, ${\mathcal A}$ is an $n$-grid iff ${\mathcal A}$ is an $(n,[n]^2)$-grid. 
                     
 If ${\mathcal I}$ is a finite set of sets, then a set $T $ is a {\bf transversal} of ${\mathcal I}$ 
  if  $T \cap I \neq \varnothing$ for every $I \in {\mathcal I}$. If $m \leq n < \omega$ and $T \subseteq n$,
  then $T$ is a transversal of $[n]^m$ iff $|T| \geq n-m+1$.
  
 The next definition  refines Definition~2.1.

\bigskip


{\sc Definition 2.3:}  Suppose that ${\mathcal A}$ is an $n$-indexed hyperspace, $|A| = \kappa$ and 
 ${\mathcal I} \subseteq {\mathcal P}(n)$. We will say that ${\mathcal A}$ is ${\mathcal I}$-{\bf twisted} if 
 
 \smallskip

$(0)$ ${\mathcal A}$  is an $(n,\{n\})$-grid
\smallskip

\noindent and either

\smallskip

$(1)$  ${\mathcal A}$ is countable, 

\smallskip

\noindent or else

\smallskip

$(2)$ $|A| = \kappa > \aleph_0$ and there is 
  a filtration $\langle A_\alpha : \alpha < \kappa \rangle$ for 
${\mathcal A}$ such that ${\mathcal A}|A_0$ is ${\mathcal I}$-twisted and whenever $\alpha < \kappa$,
$ I $ is a transversal of ${\mathcal I}$ and 
$$
B = \{x \in A_{\alpha+1} \backslash A_\alpha : \forall i < n[i \in I  \leftrightarrow [x]_i \cap A_\alpha = \varnothing   ]\},
$$ 
then $({\mathcal A}|B)\harp I$ is $({\mathcal I} \cap {\mathcal P}(I))$-twisted.  


\bigskip

If  ${\mathcal A}$ is an uncountable $n$-indexed hyperspace, then we will refer to a filtration for  
${\mathcal A}$ as in (2) of Definition~2.3 as an ${\mathcal I}$-{\bf twisted} filtration.

The following lemma relates Definitions~2.1 and~2.3.

\bigskip


{\sc Theorem 2.4}: {\em Suppose that ${\mathcal A}$ is an $(n, \mathcal I)$-grid. Then, ${\mathcal A}$ is 
$\mathcal I$-twisted iff it is twisted.}


\bigskip

{\it Proof.}  If ${\mathcal A}$ is not an $(n,\{n\})$-grid, then ${\mathcal I} = \varnothing$. Then, ${\mathcal A}$ is  $\varnothing$-twisted iff ${\mathcal A}$ is finite iff ${\mathcal A}$   has an acceptable coloring. 
Thus, it is  safe to assume that $n \in {\mathcal I}$. 

If ${\mathcal A}$ is countable, then it is both twisted and  ${\mathcal I}$-twisted. So, assume that $|A| = \kappa > \aleph_0$  
 and suppose, as an inductive hypothesis, that the theorem is valid for all smaller indexed hyperspaces. 

$(\Longleftarrow)$: Trivial.

$(\Longrightarrow)$:  Let 
$\langle A_\alpha : \alpha < \kappa \rangle$ be an ${\mathcal I}$-twisted filtration for ${\mathcal A}$. We will show that this same filtration is a twisted filtration for ${\mathcal A}$. Thus, we want to show that whenever $\alpha < \kappa$, $\varnothing \neq I \subseteq n$, $B$ is defined as in Definition 2.1(2) and ${\mathcal B} = ({\mathcal A}|B)\harp I$, then ${\mathcal B}$ is twisted. 
There are two cases.

\smallskip

{\em $I$ is not a transversal}: We claim that $B = \varnothing$. Suppose not, and let  $x \in B$. Since $I$ is not a transversal, we can pick $J \in {\mathcal I}$ such that $J \cap I = \varnothing$. Let 
$X = \bigcap\{[x]_j : j \in J\}$. Clearly, $x \in X$. Also,   $X$ is finite since $J \in {\mathcal I}$. For each $j \in J$, let $y_j \in [x]_j \cap A_\alpha$. Then, $X = \bigcap \{[y]_j : j \in J\}$. By elementarity and the finiteness of $X$, we have that $X \subseteq A_\alpha$, so that $x \in A_\alpha$, a contradiction. 

\smallskip

{\em $I$ is a transversal}: Clearly, ${\mathcal B}$ is an $(I,{\mathcal I} \cap {\mathcal P}(I))$-grid. 
Since  ${\mathcal B}$ is $({\mathcal I} \cap {\mathcal P}(I))$-twisted, then, by the inductive hypothesis, it is twisted.  \qed

\bigskip


{\sc Corollary 2.5:} {\em If ${\mathcal A}$ is an $(n, \mathcal I)$-grid, then ${\mathcal A}$ is 
$\mathcal I$-twisted iff it has an acceptable coloring.} \qed


\bigskip


{\sc Corollary 2.6:} {\em If ${\mathcal A}$ is an $n$-grid, then ${\mathcal A}$ is 
$[n]^2$-twisted iff it has an acceptable coloring.}  \qed


\bigskip

One reason for introducing Definition~2.3 and Theorem~2.4 is to be able to state the next corollary, whose main appeal is  
a characterization of the twisted $n$-grids  more resembling de la Vega's definition.

\bigskip


{\sc Corollary $2.7$}: {\em If ${\mathcal A}$ is an uncountable 
$n$-grid, then ${\mathcal A}$ is twisted   iff it has 
 a filtration $\langle A_\alpha : \alpha < \kappa \rangle$  such that ${\mathcal A}|A_0$ is twisted and whenever $\alpha < \kappa$,
$k < n$ and 
$$
B = \{x \in A_{\alpha+1} \backslash A_\alpha : [x]_k \cap A_\alpha \neq \varnothing]\},
$$ 
then $({\mathcal A}|B)\harp (n \backslash \{k\})$ is twisted.} \qed

\bigskip


\bigskip


{\bf \S3.\@ Simms's Theorems.} In \cite{sim97} Simms considered $n$-indexed hyperspaces, but allowed 
the possibility that $n$ is infinite. He   also considered some generalizations of acceptable colorings 
for these types of $n$-indexed hyperspaces.   When referring in this section to a result from \cite{sim97}, we will always be concerned just with that part of it that fits into the context of this paper. 

Suppose that   ${\mathcal I}$ is a finite set of finite sets. 
We define $\delta({\mathcal I})$, the {\bf depth} of ${\mathcal I}$, to be the  least $d < \omega$ for which there are transversals $T_0,T_1, \ldots, T_{d-1}$ of ${\mathcal I}$ such that 
$T_0 \cap T_1 \cap \cdots \cap T_{d-1} = \varnothing$.{\footnote{The term {\it depth} is borrowed from \cite[Def.~3.1]{sim97} to which it is somehow obliquely related. See Definition~3.12(a).}}
 If there are no such transversals or, equivalently, 
if there is $I \in {\mathcal I}$ such that $|I| \leq 1$, then let $\delta({\mathcal I}) = \infty$.{\footnote{ We adopt the usual conventions concerning $\infty$; for example,  $\infty -1 = \infty$ and 
$\alpha < \infty$ for every ordinal $\alpha$.}}  Some examples are: $\delta(\varnothing) = 1$; 
if $\varnothing \neq {\mathcal I} \subseteq {\mathcal P}(n)$ and $\delta({\mathcal I}) < \infty$, then $2 \leq \delta({\mathcal I}) \leq n$; $\delta([n]^2) = n$; and more generally, 
if $2 \leq m \leq n+1$, then $\delta([n]^m) = \lceil n/(m-1) \rceil$.

\bigskip

Our first goal in this section is Theorem~3.2, which extends Theorem~0.2 since \mbox{$\delta([n]^2) = n$} and also extends Lemma~1.2 since $\delta(\{n\}) = 2$ (as long as $n \geq 2$).  We give a quick proof  of Theorem~3.2 using Corollary~2.5. But first, we prove  a very simple lemma.

\bigskip


    \bigskip

   {\sc Lemma 3.1}: {\em Suppose that ${\mathcal I} \subseteq {\mathcal P}(n)$.
   If $J \subseteq n$ is a transversal of ${\mathcal I}$,
   then $\delta({\mathcal I} \cap {\mathcal P}({J})) \geq \delta({\mathcal I}) - 1$. }
   
   \bigskip
   
   {\it Proof}. If $\delta({\mathcal I} \cap {\mathcal P}(J)) = \infty$, then the conclusion is trivial, so assume that $\delta({\mathcal I} \cap {\mathcal P}(J)) \leq n$. 
   
   Suppose that ${\mathcal T} \subseteq {\mathcal P}(J)$ is a set of transversals of ${\mathcal I} \cap {\mathcal P}(J)$ such that 
   $\bigcap{\mathcal T} = \varnothing$. Let ${\mathcal T}' = \{ T \cup (n \backslash J) : T \in {\mathcal T}\} \cup \{J\}$.  
   It is easily checked that ${\mathcal T}'$ is a set of transversals of ${\mathcal I}$ and that 
   $\bigcap{\mathcal T}' = \varnothing$. To finish the proof, observe  that $|{\mathcal T}'| \leq |{\mathcal T}| + 1$. \qed

\bigskip


{\sc Theorem 3.2}:  {\em Suppose that ${\mathcal A}$ is an $(n,{\mathcal I})$-grid, $d = \delta({\mathcal I})$ and $|A| < \aleph_{d-1}$. Then ${\mathcal A}$ has an acceptable coloring.}

\bigskip

{\it Proof}. First, suppose that $d = \infty$, so there is $I \in {\mathcal I}$ such that $|I| \leq 1$.
If $I = \varnothing$, then $A$ is finite so any coloring of ${\mathcal A}$ is acceptable. If $I = \{i\}$,
then $[a]_i$ is finite for every $a \in A$, so the coloring that is constantly $i$ is acceptable. 

Next, suppose that $d < \infty$. 
We give a proof by induction on $d$. 

For the basis step,  assume that 
$d = 1$. Then ${\mathcal A}$ is finite so any coloring is acceptable..

For the inductive step,  suppose that $d \geq 2$ and that the Theorem holds for all smaller values of $d$. 

We prove by induction on the cardinal $\kappa$ that if ${\mathcal A}$ is an $(n, {\mathcal I})$-grid, 
$\delta({\mathcal I}) = d$ and $|A| = \kappa < \aleph_{d-1}$, then ${\mathcal A}$ has an acceptable coloring. 

If $\kappa \leq \aleph_0$, then Corollary~1.3 yields that ${\mathcal A}$ has an acceptable coloring.
Thus, assume that $\kappa > \aleph_0$ and that we know the result for all smaller cardinals. 

By Corollary~2.5, it suffices to show that ${\mathcal A}$ is ${\mathcal I}$-twisted. 
Let $\langle A_\alpha : \alpha < \kappa \rangle$ be any filtration for ${\mathcal A}$. 
We will show that it is ${\mathcal I}$-twisted. Let $I$ be a transversal of ${\mathcal I}$ and let $B$ be as in Definition 2.3(2).  Then $({\mathcal A}|B)\harp I$ is an $(I, {\mathcal I} \cap {\mathcal P}(I))$-grid, 
and, according to Lemma~3.1,  $\delta({\mathcal I} \cap {\mathcal P}(I)) \geq d-1$. 
Thus, by the inductive hypothesis,  $({\mathcal A}|B)\harp I$ has an acceptable coloring 
and, therefore, by Corollary~2.5, is $({\mathcal I} \cap {\mathcal P}(I))$-twisted. Hence, ${\mathcal A}$ is 
${\mathcal I}$-twisted.
  \qed
  
   \bigskip
  
  The $n$-indexed hyperspace $\big(\RR^m;E(c_0), E(c_1), \ldots, E(c_{n-1})\big)$ from the preamble 
  is an $(n, [n]^m)$-grid whenever $c_0,c_1, \ldots, c_{n-1} \in \RR^m$ are in general position.
  Thus, Theorem~3.2 implies the observation from \cite{grid} that  $\RR^m$ can be covered by $(d+1)(m-1) + 1$ sprays when $2^{\aleph_0} \leq \aleph_d$.

  Being an $(n,{\mathcal I})$-grid is  a global property of an $n$-indexed hyperspace. This can modified into a more local property as follows. Let ${\mathcal A}$ be an $n$-indexed hyperspace. 
  For each $a \in A$, let ${\mathcal I}(a) = \{I \subseteq n : |\bigcap_{i \in I}[a]_i| < \aleph_0\}$, and then let ${\mathcal I}({\mathcal A}) = \bigcap\{{\mathcal I}(a) : a \in A\}$. Thus, ${\mathcal I}({\mathcal A})$ is the set of all those $I \subseteq n$ such that ${\mathcal A}$ is an $(n,\{I\})$-grid. By Theorem~3.2,
  if  ${\mathcal A}$ is an $n$-indexed hyperspace, $d = \delta({\mathcal I}({\mathcal A}))$ and $|A| < \aleph_{d-1}$, then ${\mathcal A}$ has an acceptable coloring. Theorem~3.2 implies Corollary~3.3, which is a local version of Theorem~3.2. 
   Corollary~3.3 is   slightly stronger than  Simms's theorem 
 \cite[Theorem~3.2]{sim97}.  The relation between Corollary~3.3 and Simms's theorem is clarified at the end of this section.
\bigskip

  \bigskip
  
  {\sc Corollary 3.3:} {\em Suppose that ${\mathcal A}$ is an $n$-indexed hyperspace,  $1 \leq d \leq  \delta({\mathcal I}(a))$ for each $a \in A$, and  $|A| < \aleph_{d-1}$. Then ${\mathcal A}$ has an acceptable coloring.}
  
  \bigskip
  
  {\it Proof.} Let ${\mathcal I}_0, {\mathcal I}_1, \ldots, {\mathcal I}_m$ be all those ${\mathcal I} \subseteq {\mathcal P}(n)$ for which $\delta({\mathcal I}) \geq d$. For each $j \leq m$, let 
  ${\mathcal A}_j = \{a \in A : {\mathcal I}(a) = {\mathcal I}_j\}$. Then $A_0,A_1, \ldots, A_m$ 
  partitions $A$ (but with the possibility that some $A_j = \varnothing$).
  Clearly,  ${\mathcal A} | A_j$ is an $(n, {\mathcal I}_j)$-grid, so, by Theorem~3.2,
  ${\mathcal A}|A_j$ has an acceptable coloring $\varphi_j$. Then, $\varphi = \bigcup_j\varphi_j$ is an acceptable coloring of ${\mathcal A}$. \qed

\bigskip

The concept of an $n$-cube  will be generalized. Suppose that $A = A_0  \times A_1 \times \cdots \times A_{m-1}$,
where $A_0,A_1, \ldots, A_{m-1}$ are arbitrary nonempty sets. If $S \subseteq m < \omega$, then $S$ {\bf induces} the 
equivalence relation $E$ on $A$, where $E$ is  such that if $x,y \in A$, then $\langle x,y \rangle \in E$ iff 
 $x_j = y_j$ whenever $j \in m \backslash S$. If $m < \omega$ and $\vS = \langle S_0,S_1, \ldots, S_{n-1} \rangle$ is an $n$-tuple of subsets of $m$,  
 then the $\vS$-{\bf cube for} $A$ is the $n$-indexed hyperspace  ${\mathcal A} = (A;E_0,E_1, \ldots, E_{n-1})$, where each $E_i$ is induced by $S_i$. The  ${\vS}$-{\bf cube over} $X$ is the $\vS$-cube for $X^m$. 
 Observe that the $n$-cube $X^n$ is exactly the ${\vS}$-cube over $X$, where 
 $\vS = \langle \{0\}, \{1\}, \ldots, \{n-1\} \rangle$. If $I$ and $M$ are   finite sets 
 and $\vS = \langle S_i : i \in I \rangle$ is an $I$-tuple of subsets of $M$, then 
 the notions of an $\vS$-cube  and an $\vS$-cube over $X$ have the obvious definitions.
  Also, for such an $\vS$, if $J \subseteq I$, then $\vS \harp J = \langle S_i : i \in J \rangle$.

   In these definitions when we have an $n$-tuple $\vS$ of subsets of $m$,  it will always be understood what $m$ is, and we leave it implicit.

    If  ${\mathcal I}$ is a finite set of   sets, then 
   we define the {\bf transversal number} of ${\mathcal I}$, and denote it by 
$\cov({\mathcal I})$, to be the least  
 cardinality of a transversal of ${\mathcal I}$. If $\varnothing \in {\mathcal I}$, then ${\mathcal I}$ does not have a transversal, so we conventionally let $\cov({\mathcal I}) = \infty$.  Note that $\tau({\mathcal I}) = 0$ iff ${\mathcal I} = \varnothing$. If $\vS = \langle S_0,S_1, \ldots, S_{n-1} \rangle$ is an $n$-tuple of  sets, then a {\bf transversal} of $\vS$ is a transversal of 
 $\{S_0,S_1, \ldots, S_{n-1} \}$ and we let 
 $\cov(\vS) = \cov(\{S_0,S_1, \ldots, S_{n-1}\})$. If $\vS$ is an $n$-tuple of nonempty subsets of $m$, then $\tau(\vS) \leq \min(m,n)$. 
   
   If $\vS = \langle S_0,S_1, \ldots, S_{n-1} \rangle$ is an $n$-tuple of finite sets, then we let 
   ${\mathcal I}(\vS) = \{I \subseteq n : \bigcap_{i \in I}S_i = \varnothing\}$. The point of this definition is that every $\vS$-cube  is an $(n, {\mathcal I}(\vS))$-grid and, moreover, whenever ${\mathcal A}$ is an 
   $\vS$-cube over an infinite set, $a \in A$, $I \subseteq n$, then $\bigcap_{i \in I}[a]_i$ is finite iff $I \in {\mathcal I}(\vS)$. The next  lemma 
   describes the relationship between  the transversal number of 
   $\vS$ and the depth of ${\mathcal I}(\vS)$.
   
   \bigskip
   
   
   {\sc Lemma 3.4}: {\em  If $\vS$ is an $n$-tuple of finite sets, then
   $\tau(\vS) = \delta({\mathcal I}(\vS))$.}
   
   \bigskip
   
   {\it Proof}. Let $\vS = \langle S_0,S_1, \ldots,S_{n-1} \rangle$, $t = \tau(\vS)$ and $d = \delta({\mathcal I}(\vS))$. 
   
   First, notice that $t = \infty$ iff some $S_i = \varnothing$ iff some $\{i\} \in {\mathcal I}(\vS)$ iff $\delta({\mathcal I}(\vS)) = \infty$. Hence, we assume that $d,t < \infty$. 
   
   Let $\{a_0,a_1, \ldots,a_{t-1}\}$ be a transversal of $\vS$. For each $k<t$, 
   let $I_k = \{i < n : a_k \not\in S_i\}$. Then, each $I_k$ is a transversal of ${\mathcal I}(\vS)$ and 
   $\bigcap_{k < t}I_k = \varnothing$, so that $d \leq t$.
   
   Conversely, let $\{T_0,T_1, \ldots, T_{d-1}\}$ be a set of transversals of  ${\mathcal I}(\vS)$ such that   $T_0 \cap T_1 \cap \cdots \cap T_{d-1} = \varnothing$. 
   For each $j < d$, let $s_j \in \bigcap_{i \in T_j}S_i$. Then $\{s_j : j < d\}$ is a transversal 
   of $\vS$, so that $t \leq d$. \qed
   
     \bigskip
   
   
   {\sc Definition 3.5:}    Suppose that ${\mathcal A}$ and ${\mathcal B}$ are, respectively, $n$-indexed and $d$-indexed hyperspaces and that  $\beta : n \into d$.  We say that a function $f : B \into A$ is a $\beta$-{\bf parbedding} of ${\mathcal B}$ into ${\mathcal A}$ if it is one-to-one and 
   for   $x,y \in B$  and $i < n$, 
    \begin{equation} \tag{$*$}
    [x]_{\beta(i}) = [y]_{\beta(i)} \Longrightarrow [f(x)]_i = [f(y)]_i.
    \end{equation}
   If   $f$ is a $\beta$-parbedding of ${\mathcal B}$ into ${\mathcal A}$ for some $ \beta$, then $f$ is a {\bf parbedding} of ${\mathcal B}$ into ${\mathcal A}$, in which case   we say that ${\mathcal B}$ is {\bf parbeddable} into ${\mathcal A}$ or that ${\mathcal A}$ parbeds ${\mathcal B}$.
   
   \bigskip
   
    Note that every embedding is a $\beta$-parbedding, where  $\beta$ is the identity function. 
    
    parbeddability is transitive. In fact, if ${\mathcal A}_0, {\mathcal A}_1, {\mathcal A}_2$ are, respectively,
    $n_0$-, $n_1$-, $n_2$-indexed hyperspaces, $\alpha : n_1 \into n_0$, $\beta : n_2 \into n_1$ 
    and $f : A_0 \into A_1$ and $g : A_1 \into A_2$ are, respectively, an $\alpha$-parbedding of ${\mathcal A}_0$ into ${\mathcal A}_1$ and a $\beta$-parbedding of ${\mathcal A}_1$ into ${\mathcal A}_2$, 
    then $gf$ is a $\beta\alpha$-parbedding of ${\mathcal A}_0$ into ${\mathcal A}_2$. 
    
    If $X$ is infinite and $2 \leq n < \omega$, then the $(n+1)$-cube over $X$ is $\beta$-parbeddable into the $n$-cube over $X$, where $\beta : n \into n+1$ is the identity function.

       \bigskip
       

 {\sc Lemma 3.6:} {\em Suppose that ${\mathcal A}$ and ${\mathcal B}$ are, respectively, $n$-indexed and $d$-indexed hyperspaces and that $\mathcal A$ parbeds ${\mathcal B}$.   If  ${\mathcal A}$ has 
   an acceptable coloring, then  so does ${\mathcal B}$.}
   
   \bigskip
   
   {\it Proof}.  Suppose that ${\mathcal A}$ and ${\mathcal B}$ are, respectively, $n$-indexed and $d$-indexed hyperspaces. Let $f : B \into A$ be an $\beta$-parbedding of ${\mathcal B}$ into ${\mathcal A}$. Let $\chi : A \into n$ be an acceptable coloring for ${\mathcal A}$, and then let  $\psi = \beta \chi f$. Clearly, $\psi : B \into d$, so $\psi$ is a coloring for ${\mathcal B}$. 
We claim that $\psi$ is  acceptable.    

   For a  contradiction, suppose that $\psi$ is not acceptable. Thus, we have  $b \in B$, $j < d$ and an infinite $X \subseteq \{x \in [b]_j : \psi(x) = j\}$. Let  $a = f(b)$ and $Y = f[X]$.   By the Pigeon Hole Principle, we can assume that $i < n$ is such that $\chi$ is constantly $i$ on $Y$.     Thus, $\beta(i) = j$ so that 
  $[x]_{\beta(i)} = [b]_{\beta(i)}$ for all $x \in X$. But then $[y]_i = [a]_i$ for all $y \in Y$ by $(*)$ of Definition~3.5. Since $Y$ is infinite and $\chi$ is acceptable, this is impossible. \qed

   \bigskip

 The next lemma gives some of the main examples of parbeddability.
 
 \bigskip
 

 {\sc Lemma 3.7}: {\em Suppose that  $1 \leq m,n  < \omega$,  ${\vS}$ is an $n$-tuple of
 nonempty  subsets  of $m$, 
$d = \cov({\vS}) \geq 1$ and $X$ is any set. Then the $d$-cube over $X$ is parbeddable into the $\vS$-cube over $X$.} 

 \bigskip

  {\it Proof}. Let $m,n, \vS, d$ and $X$ be as given.
    Let ${\mathcal A}$ be the $\vS$-cube over $X$.  
  Let $T = \{t_0,t_1, \ldots, t_{d-1}\} \subseteq m$ be a transversal of $\vS$. Let $\beta : n \into d$ be such that $t_{\beta(i)} \in S_i$ for each $i < n$. There is such a $\beta$ since $T$ is a transversal. Fix $c \in X$. Define $f : X^d \into X^m$ so that if $x \in X^d$,  then 
    \begin{displaymath}
  f(x)_k = \left\{
  \begin{array}{ll} x_j & \textrm{if~} k = t_j \\ c & \textrm{otherwise.} 
  \end{array} \right.
  \end{displaymath}
  We will show that $f$ is a $\beta$-parbedding of $X^d$ into ${\mathcal A}$ by proving that if $x,y \in X^d$ and $i < n$, then $(*)$ of Definition~3.5 holds. 
  \begin{eqnarray}  
[x]_{\beta(i)} = [y]_{\beta(i)} & \Longrightarrow & \forall j < d \big(j \neq \beta(i) \into x_j = y_j \big) \nonumber \\
& \Longrightarrow & \forall k<m \big(k \neq t_{\beta(i)} \into f(x)_k = f(y)_k\big) \nonumber \\
& \Longrightarrow & \forall k < m \big(k \not\in S_i \into    f(x)_k = f(y)_k\big) \nonumber \\
& \Longrightarrow & [f(x)]_i = [f(y)]_i \ . \nonumber 
\end{eqnarray}  

  \vspace{-19pt} 
  \qed
  \vspace{19pt}

   
   \bigskip

The next theorem  implies  Kuratowski's Theorem 0.1  because 
$\tau(\langle \{0\},$ $ \{1\}, \ldots, \{n-1\} \rangle ) = n$.{\footnote{There is more to this story. The first  three papers published in \cite{fm} contain a sequence of three successively stronger theorems. 
In the first one, Sierpi\'{n}ski \cite{sie51} proves Theorem~0.1 with $n=3$; in the second,
Kuratowski \cite{kur51} proves his Theorem~0.1; and in the third one, Sikorski \cite{sik51} 
proves the special case of Theorem~3.8 in which $1 \leq k \leq m$, $n = \binom{m}{k}$ and $\vS$ is an $n$-tuple of all the $k$-element subsets of $m$. Theorem~0.1 with $n=3$ was proved twice more by Sierpi\'{n}ski \cite{sie52}, \cite{sie51a}. The special case of Theorem~0.1 with $n=2$ was also proved by Sierpi\'nski \cite{sie51b}. More about these historical developments can be found in \cite{ejm94}, \cite{grid} and especially \cite{sim91}.}} Theorem~3.8 is a consequence of the somewhat arcane Theorem~4.3 of \cite{sim97}. Erd\H{o}s, Jackson and Mauldin in \cite[Coro.\@ 7]{ejm94} made a further generalization of Theorem~3.8 that allowed for structures even more general than 
$n$-indexed hyperspaces.{\footnote{The reference in \cite{ejm94} that is identified there by [Sm2] is apparently  a preliminary version of \cite{sim97}.}}

 \bigskip
 
 
 {\sc Theorem 3.8:}  (Simms \cite{sim97})  {\em Suppose that  $1 \leq m,n  < \omega$,  ${\vS}$ is an $n$-tuple of nonempty 
  subsets  of $m$, 
$d = \cov({\vS}) \geq 1$, and $X$ is a set.  Then the 
 ${\vS}$-cube over $X$ has an acceptable coloring iff $|X| < \aleph_{d-1}$.}
 

   \bigskip
   
   {\it Proof}.    Let $m,n,\vS,d$  be as given, and let ${\mathcal A}$ be the $\vS$-cube over $X$. 
   
      \smallskip
   
   $(\Longrightarrow)$:    Suppose that $|X| \geq \aleph_{d-1}$. By Lemma~3.7, the $d$-cube $X^d$ is parbeddable into ${\mathcal A}$. Theorem~0.1 implies that the $d$-cube over $X$ does not have an acceptable coloring and therefore, by Lemma~3.6, neither does ${\mathcal A}$.
   
   \smallskip
   
   $(\Longleftarrow)$: Suppose that $|X| < \aleph_{d-1}$. As already noted, ${\mathcal A}$ is an $(n, {\mathcal I}(\vS))$-grid. By Lemma~3.4, $d = \tau(\vS) = \delta({\mathcal I}(\vS))$, so that ${\mathcal A}$ has an acceptable coloring by Theorem~3.2.  \qed
   
   \bigskip
   
   As an example, consider the $n$-tuple $\vS = \langle S_0,S_1, \ldots, S_{n-1} \rangle$, where $S_i = n \backslash \{i\}$. Then $\tau(\vS) = 2$, so that the $\vS$-cube over $\RR$ does  not have an acceptable coloring. 
   The $\vS$-cube over $\RR$ is embeddable into $(\RR^n;E(c_0),E(c_1), \ldots, E(c_{n-1})$ from the preamble, implying that $\RR^n$ cannot be covered by $n$ sprays.
   
   We prove one more result along these lines.
   
   \bigskip   
   
     {\sc Theorem 3.9:}  {\em Suppose that  $1 \leq k < m  < \omega$ and $n = \binom{m}{k}$. Let  ${\vS}$ be an $n$-tuple  of all $k$-element subsets of $m$. Let 
${\mathcal A}$ be an $\vS$-cube, where  $A = X_0 \times X_1 \times \cdots \times X_{m-1}$.   Then ${\mathcal A}$ has an acceptable coloring iff there is $d < m-k+1$ such that 
$|\{j < m : |X_j| < \aleph_d\}| \geq d+k$.}

   \bigskip

   {\it Proof.} $(\Longleftarrow)$: Let $d$  be as in the Theorem. Let 
   $J = \{j < m : |X_j| < \aleph_d\}$ and $I = \{i < n : S_i \subseteq J\}$. Let 
   ${\mathcal B}$  be the $\vS \harp I$-cube, where  
   $B = \prod_{j \in J}X_j = \{x \harp J : x \in A\}$. 
   Then,   $\tau(\vS \harp I) = |J| - k +1 \geq d +1$ and $|X_j| < \aleph_d$ for each $j \in J$. Theorem~3.8 implies that ${\mathcal B}$ has an acceptable coloring $\psi : B \into I$. Let $\chi :  A \into I$ be such that 
   for $x \in A$, $\chi(x) = \psi(x \harp J)$. Then, $\chi$ is an  
 acceptable coloring of ${\mathcal A} \harp I$ and, therefore,  is also 
 an acceptable coloring of ${\mathcal A}$.

\smallskip  

$(\Longrightarrow)$: The special case when $m = n$ and $k = 1$ is known.{\footnote{See \cite[Prop.\@ 2.151]{sim91}. As mentioned in \cite{sim91}, there was some confusion about the attribution.
However, the reference given in \cite{sim91} does not clarify it  since it does not correspond to an entry 
in the References of \cite{sim91}.}} 
 Thus, the $n$-cube $\aleph_0 \times \aleph_1 \times \cdots \times \aleph_{n-1}$ does not have an acceptable coloring. 

Suppose that whenever $d < m-k+1$,  
then $|\{j < m : |X_j| < \aleph_d\}| \leq d+k - 1$. Since $\tau(\vS) = m-k+1$, it has a transversal 
$T = \{t_0,t_1, \ldots, t_{m-k}\}$. As in the proof of Lemma~3.7, $X_{t_0} \times X_{t_1} \times \cdots \times X_{t_{m-k}}$ is parbeddable into ${\mathcal A}$. Assuming, without loss of generality, 
that $|X_{t_0}| \leq  |X_{t_1}| \leq \ \cdots \leq |X_{t_{m-k}}|$, we then have that $|X_{t_j}| \geq \aleph_{t_j}$ for $j \leq m-k$, so that the $(m-k+1)$-cube $\aleph_0 \times \aleph_1 \times \cdots \times 
\aleph_{m-k}$ is parbeddable into ${\mathcal A}$. Thus, by Lemma~3.6 and the result mentioned in the previous paragraph, we have that ${\mathcal A}$ does not have an acceptable coloring. 
\qed

\bigskip

{\it Question 3.10}: Is there a generalization of Theorem~3.9 that applies to all $\vS$-cubes?

      \bigskip
   
Theorem~3.8 is primarily about infinite $\vS$-cubes. Nevertheless, 
 finite $\vS$-cubes will play a significant role in  the next section. 
  \bigskip
  

   \bigskip

  
  As already mentioned, Corollary~3.3 is somewhat stronger than \cite[Theorem~3.2]{sim97}. 
  The remainder of this  section is devoted to discussing the relation between  these results. 
  
  Part (a) of the following  definition is taken directly from \cite[Def.\@ 3.1]{sim97}, and this definition naturally suggests the one in (b).

\bigskip

{\sc Definition 3.12:} (\cite[Def.\@ 3.1]{sim97}) Suppose that $m < \omega$, $d \leq n < \omega$ and ${\mathcal A}$ is  An $n$-indexed hyperspace.

(a)  ${\mathcal A}$ is 
$m$-{\bf fine to depth} $d$ if whenever $a \in A$ and $\pi$ is a permutation of $n$, then there 
are $0 = i_0 \leq i_1 \leq i_2 \leq \cdots \leq i_{d} < n$ such that  whenever 
$k < d$,  then $|\bigcap\{[a]_{\pi(j)} : i_k \leq j  \leq i_{k+1}\}| \leq m$. 

(b) ${\mathcal A}$ is 
{\bf fine to depth} $d$ if whenever $a \in A$ and $\pi$ is a permutation of $n$, then there 
are $0 = i_0 \leq i_1 \leq i_2 \leq \cdots \leq i_{d} < n$ such that  whenever 
$k < d$,  then $\bigcap\{[a]_{\pi(j)} : i_k \leq j  \leq i_{k+1}\}$ is finite.

\bigskip

Observe that if ${\mathcal A}$ is $m$-fine to depth $d$, $m \leq m' < \omega$ and $d' \leq d$, 
then ${\mathcal A}$ is $m'$-fine to depth $d'$ and also is fine to depth $d$.    In stating \cite[Theorem~3.2] {sim97}, Simms does not use the notion of depth that was used in our Theorem~3.2 and Corollary~3.3, but uses instead the notion defined in Definition~3.12(a). If ${\mathcal I} \subseteq 
{\mathcal P}(n)$, then we will say that ${\mathcal I}$ is {\bf dandy to depth} $d$ if, 
for every permutation $\pi $  of $n$,  there 
are $0 = i_0 \leq i_1 \leq i_2 \leq \cdots \leq i_{d} < n$ such that  for every 
$k < d$ there is $I \in {\mathcal I}$ such that  $I \subseteq  \{\pi(j) : i_k \leq j \leq   i_{k+1}\}$.

\bigskip


{\sc Lemma 3.13}: {\em Suppose that ${\mathcal I} \subseteq {\mathcal P}(n)$ and $d < \omega$.
Then $d < \delta({\mathcal I})$ iff  ${\mathcal I}$ is dandy to depth $d$.} 

\bigskip

{\it Proof.} First, suppose that $\delta({\mathcal I})= \infty$. Thus, there is $I \in {\mathcal I}$ such that $|I| \leq 1$. 
Let $d < \omega$ and  $\pi$ be a permutation  of $n$. If $\varnothing \in {\mathcal I}$, then let $0 = i_0 = i_1 = \cdots = i_d$. Otherwise, let $j < n$ be such that $\{\pi(j)\} \in {\mathcal I}$ and let 
$0 = i_0 \leq i_1 = i_2 = \cdots = i_d = j$. Either way, we see that ${\mathcal I}$ is dandy to depth $d$.

Next, if $\delta({\mathcal I}) = 1$, then ${\mathcal I} = \varnothing$, so that we easily see that 
${\mathcal I}$ is dandy to depth $0$ and not dandy to depth $1$. 
So, assume that $2 \leq \delta({\mathcal I}) \leq n$.

\smallskip

$(\Longrightarrow)$: Suppose that $d = \delta({\mathcal I}) -1$. Let $\pi$ be a permutation of $n$. Without loss of generality,  assume that $\pi$ is the identity permutation. Define the sequence 
 $0 = i_0 < i_1 < i_2 < \cdots < i_{e} < n$  such that $e$ is as large as possible and whenever 
$k < e$, then $i_{k+1}$ is the least for which there is $I \in {\mathcal I}$ such that 
$I \subseteq  \{j < n: i_k \leq j \leq i_{k+1}\}$. Thus, ${\mathcal I}$ is dandy to depth $e$. For $k < e$, let $I_k$ be the interval $[i_k,i_{k+1})$ 
and let $I_{e} = [i_e,n)$. Then let $T_k = n \backslash I_k$ for $k \leq e$. Each $T_k$ is a transversal of ${\mathcal I}$ and $T_0 \cap T_1 \cap \cdots \cap T_{e} = \varnothing$. Thus, $e \geq d$, proving that ${\mathcal I}$ is dandy to depth $d$.

\smallskip
$(\Longleftarrow)$: We wish to show that if $d = \delta({\mathcal I})$, then ${\mathcal I}$ is not 
dandy to depth $d$. The proof is by induction on $d$. 

$d = 2$: For a contradiction, assume that ${\mathcal I}$ is dandy to depth $2$. Let $T_0, T_1$ be transversals of ${\mathcal I}$ such that $T_0 \cap T_1 = \varnothing$.
We can assume that $T_0 \cup T_1 = n$. 
Let $\pi$ be a permutation of $n$ such that if $i \in T_0$ and $j \in T_1$, then $\pi(i) < \pi(j)$. 
Without loss of generality, assume that $\pi$ is the identity permutation. 
Let 
 $0 = i_0 < i_1 < i_2 < n$  demonstrate that ${\mathcal I}$ is dandy to depth $2$; that is,   there are $I_0,I_1 \in {\mathcal I}$ such that 
$I_0 \subseteq  [0,i_1]$ and $I_1 \subseteq [i_1,i_2]$. Then $i_1 \in T_1)$ so that $T_0 \cap I_0 
\subseteq T_0 \cap [i_1,i_2] = \varnothing$, contradicting that $T_0$ is a transversal.  

\smallskip

For the inductive step, let $2 < d \leq n$ and assume that for all smaller $d$ we have the result. 
For a contradiction, assume that ${\mathcal I}$ is dandy to depth $d$.
Let $T_0,T_1, \ldots, T_{d-1}$ be transversals of 
${\mathcal I}$ such that $T_0 \cap T_1 \cap \cdots \cap T_{d-1} = \varnothing$. We can assume that 
$T_0 \cup T_1 \cup \cdots \cup T_{d-1} = n$. 
Let $\pi$ be a permutation of $n$ such that whenever $k < d$, $i \in T_k$ and $j \not\in T_0 \cup T_1 \cup \cdots \cup T_k$, then 
$\pi(i) < \pi(j)$. Without loss of generality, assume that $\pi$ is the identity permutation. 
Let 
 $0 = i_0 < i_1 < \cdots < i_d < n$  demonstrate that ${\mathcal I}$ is dandy to depth $d$; that is,   there are $I_0,I_1, \ldots, I_{d-1} \in {\mathcal I}$ such that 
$I_k \subseteq   [i_k,i_{k+1}]$ for $k < d$. Then, $i_1 \in T_d$, so that $i_1,i_2, \ldots, i_{d-1} \in T_d$.
Thus, $i_1 < i_2 < \cdots < i_{d-1}$ demonstrate that ${\mathcal I} \cap {\mathcal P}([i_1,n))$ 
is dandy to depth $d-1$. This implies that ${\mathcal I} \cap {\mathcal P}(T_d)$ is dandy to depth 
$d-1$. Then, by the inductive hypothesis, $\delta({\mathcal I} \cap {\mathcal P}(T_d)) \geq d$. However,
$T_0 \cap T_d, T_1 \cap T_d, \ldots, T_{d-1} \cap T_d$ are $d-1$ transversals of ${\mathcal I} \cap {\mathcal P}(T_d)$
 whose intersection is $\varnothing$, thereby showing 
the contradiction that  $\delta({\mathcal I} \cap {\mathcal P}(T_d)) \leq d-1$. \qed

\bigskip

{\sc Corollary 3.14:} ({\it cf.\@} \cite[Theorem~3.2]{sim97}) {\em Suppose that  ${\mathcal A}$ is an $n$-indexed hyperspace that is fine to depth $d$. If  $|A| < \aleph_d$, 
then ${\mathcal A}$ has an acceptable coloring.}

\bigskip

{\it Proof}. It follows from Lemma~3.13 that ${\mathcal A}$ is fine to depth $d$ iff 
$\delta({\mathcal I}(a)) > d$. Hence, by Corollary~3.3, ${\mathcal A}$ has an acceptable coloring. \qed

\bigskip 

The hypothesis of the corollary is implied by the weaker one that for some $m < \omega$, ${\mathcal A}$ is $m$-fine to depth $d$. It is exactly this latter hypothesis that  Theorem~3.2 of \cite{sim97} has when it is restricted to our context. Corollary~3.14 (and its equivalent Corollary~3.3) is strictly stronger 
than \cite[Theorem~3.2]{sim97} as the following example shows. Let $E$ be an equivalence relation 
on an infinite set $A$ all of whose equivalence classes are finite and for which there are arbitrarily large finite equivalence classes. Then, the $n$-indexed hyperspace ${\mathcal A} = (A;E,E, \ldots, E)$ 
is fine to depth $n$ but for no $m < \omega$ is it $m$-fine to depth~$1$.

\bigskip


 {\bf \S4.\@ Extending de la Vega's theorem.}     As its title suggests, 
 this section's  main purpose  is to extend de la Vega's Theorem~0.3 from $n$-grids to $n$-indexed hyperspaces. This will be done in Theorem~4.2. At the same time, the hypothesis of Theorem~0.3 will be weakened, yielding Corollary~4.4. In Theorem~4.5 we give a modification of Theorem~4.2 that 
 restricts the cardinality of the indexed hyperspaces. 
 
 If ${\mathcal A}$ and ${\mathcal B}$ are $n$-indexed hyperspaces, then a {\bf weak embedding} of ${\mathcal B}$ into ${\mathcal A}$ is a one-to-one function $f : B \into A$ for which there is a permutation $\pi$ of $n$ such that whenever $x,y \in B$ and $i < n$, then 
 $$
 [x]_{\pi(i)} = [y]_{\pi(i)} \Longleftrightarrow [f(x)]_i = [f(y)]_i.
 $$
 If there is a weak embedding of ${\mathcal B}$ into ${\mathcal A}$, then we say that ${\mathcal B}$ is {\bf weakly embeddable} into ${\mathcal A}$ or that ${\mathcal A}$ {\bf weakly embeds} ${\mathcal B}$. Obviously, if ${\mathcal A}$ weakly embeds ${\mathcal B}$ and ${\mathcal A}$ has an acceptable coloring, then so does ${\mathcal B}$. Every embedding is a weak embedding. If ${\mathcal B}$ is an  
 $n$-cube over $X$, then ${\mathcal B}$ is weakly embeddable into ${\mathcal A}$ iff ${\mathcal B}$ is embeddable into ${\mathcal A}$. Every  weak embedding of ${\mathcal B}$ into ${\mathcal A}$ is  a 
 parbedding; in fact, if $\pi$ is a permutation that witnesses that $f$ is a  weak emebbeding, then $f$ is a $\pi$-parbedding. 
 
 For any linearly ordered set $X$ (for example, any $X \subseteq \omega$) and $n < \omega$, let 
 $\langle X \rangle^n$ be the set of strictly increasing $n$-tuples from $X$. 
 Define the  $n$-{\bf halfcube} to be  the $n$-grid ${\mathcal A}|\langle \omega \rangle^n$, 
 where ${\mathcal A}$ is  the $n$-cube over $\omega$. 
 If $1 \leq m, n < \omega$ and $\vS$ is an $n$-tuple of finite subsets of $m$, 
 then we define the 
 $\vS$-{\bf halfcube} to be the $n$-indexed hyperspace ${\mathcal A}|\langle \omega \rangle^m$, where ${\mathcal A}$ is  the ${\vS}$-cube over $\omega$. Thus, the $ n$-halfcube is just the $\langle \{0\}, \{1\}, \ldots, \{n-1\} \rangle$-halfcube. 
  For any $\vS$, the $\vS$-halfcube embeds all finite $\vS$-cubes; however,  there 
 are $n$-indexed hyperspaces that embed all finite $\vS$-cubes but do not embed the
 $\vS$-halfcube. 
                 
               Recall that Infinite Ramsey's Theorem asserts that whenever ${\mathcal P}$ is a finite 
               partition of $\langle \omega \rangle^n$, then there are $P \in {\mathcal P}$ and an infinite $X \subseteq \omega$ such that $\langle X \rangle^n \subseteq P$. We will need the following canonical version of Ramsey's Theorem due to Erd\H{o}s \& Rado \cite{er50}.                
               \bigskip
               
               {\sc Lemma 4.1:} (Erd\H{o}s-Rado) {\em Let $n < \omega$ and let $E$ be any equivalence relation on $\langle \omega\rangle^n$. Then there are $I \subseteq n$ and an infinite $X \subseteq \omega$  such that 
               $E \cap (\langle X \rangle^n)^2$ is the equivalence relation on $\langle X \rangle^n$ induced by $I$.}
               
               \bigskip

            {\sc Theorem 4.2}: {\em Suppose that  ${\mathcal A}$ is an $n$-indexed hyperspace that does not weakly embed any $\vS$-halfcube, where $\vS$ is an   $n$-tuple   of  nonempty subsets of $n$. 
             Then ${\mathcal A}$ has an acceptable coloring.}    
          
           \bigskip
    
    {\it Proof.} We prove the theorem by induction on the cardinality of ${\mathcal A}$. First, assume that ${\mathcal A}$ is countable.  Then $[a]_0 \cap [a]_1 \cap \cdots \cap [a]_{n-1}$ is finite 
    for every $a \in A$ as otherwise each $\langle n,n, \ldots, n \rangle$-halfcube would be embeddable into ${\mathcal A}$.  
     By Lemma~1.2, ${\mathcal A}$ has an acceptable coloring. 
   
     Next, suppose that ${\mathcal A}$ has cardinality $\kappa > \aleph_0$ and assume, as an inductive hypothesis, that the theorem is valid when restricted to $n$-indexed hyperspaces of smaller cardinality..
    We will prove that ${\mathcal A}$ is twisted, which, by Theorem~2.2, implies that ${\mathcal A}$ has an acceptable coloring.  Thus, it suffices to show that there is a twisted filtration for ${\mathcal A}$. We will prove that every filtration for ${\mathcal A}$ is twisted. 
    
     Let $\langle A_\alpha : \alpha < \kappa \rangle$ be a filtration for ${\mathcal A}$. Clearly, ${\mathcal A}|A_0$ satisfies the hypothesis of the Theorem and $|A_0| < \kappa$; hence, by the inductive hypothesis, ${\mathcal A}|A_0$ is twisted. Next,      consider $\alpha < \kappa$ and nonempty $I \subseteq n$, and then let  $B$ be as in (2) of Definition~2.1 and ${\mathcal B} = ({\mathcal A}|B)\harp I$. We wish to show that ${\mathcal B}$ is twisted or, equivalently, that      ${\mathcal B}$ has an acceptable coloring. To do so, we will use the inductive hypothesis and then prove: whenever $\vR$ is an $I$-tuple of  nonempty subsets of $I$, 
    then ${\mathcal B}$ does not weakly embed the ${\vR}$-halfcube.

     For a contradiction, suppose that $\vR$ is an $I$-tuple of  nonempty subsets of $I$ and $f : \langle \omega \rangle^I \into B$ is a weak embedding of  the ${\vR}$-halfcube 
    into ${\mathcal B}$.  For notational convenience and 
    without loss of generality, we assume that $I = m > 0$ and that $f$ is actually an embedding of  the $\vR$-halfcube into ${\mathcal B}$. It must be that $m < n$, as otherwise $I = n$ and $f$ would be an embedding of the $\vR$-halfcube into ${\mathcal A}$.

      We define  a function $g : \langle \omega \rangle^n \into A_\alpha$ by recursion. 
           
     For  $m \leq i < n$ and $c \in \langle \omega \rangle^m$,   let $a_{c,i} \in A_\alpha$ be such that $[f(c)]_i = [a_{c,i}]_i$. 
     The function $g$ will be obtained as the union of an increasing  sequence $g_0 \subseteq g_1 \subseteq g_2 \subseteq \cdots$, where, for each $r < \omega$, $g_r : \{d \in \langle \omega \rangle^n : d_{n-1} <r\} \into A_\alpha$.  There is no choice for $g_r$ when $r < n$ since  each domain is $\varnothing$.
     Now suppose we have  $g_r$ and wish to get $g_{r+1}$.

      Let 
      $$
      X =  \{g_r(d) : d_{n-1} < r\} \cup \{a_{d,i} : m \leq i < n {\mbox{ and }} d_{n-1} = r\}.
      $$ 
      Since $X$ is a finite subset of $A_\alpha$,  by elementarity, we can get $g_{r+1} \supseteq g_r$ such that whenever $d,e \in \langle \omega \rangle^n$ and $d_{n-1} = r = e_{n-1}$, then:
     \begin{itemize}
     
     \item[(1)]  $g_{r+1}(d) \in A_\alpha \backslash X$;
     
     \item[(2)] if $d \neq e$, then $g_{r+1}(d) \neq g_{r+1}(e)$;
     
     \item[(3)] if $x \in X$ and $i < n$,  then 
       $
     [g_{r+1}(d)]_i = [x]_i \Longleftrightarrow [f(d \harp m)]_i = [x]_i
     $;
     
     \item[(4)] if $i < n$, then 
      $
     [g_{r+1}(d)]_i = [g_{r+1}(e)]_i \Longleftrightarrow [f(d \harp m)]_i = [f(e \harp m)]_i$.

     \end{itemize}
     
     By $(1)$ and $(2)$, this defines a one-to-one function $g : \langle \omega \rangle^n \into A_\alpha$. We claim:
     
     \begin{quote}  \hspace{-30pt} $(*)$ \, {\em For each $i < n$ there is $j < n$ such that whenever $d,e \in \langle \omega \rangle^n$ are such that  $d_k = e_k$ whenever $j \neq k < n$, then $[g(d)]_i = [g(e)]_i$.}  \end{quote}  The proof of the claim divides into  two cases depending on whether or not $i < m$.
     
     \smallskip
     
     $i<m$: Let $j \in R_i$, which is possible since $R_i \neq \varnothing$. Thus, $j < m$. 
     Let $d,e \in \langle \omega \rangle^n$ be such that  $d_k = e_k$ whenever $j \neq k < n$, intending to prove that $[g(d)]_i = [g(e)]_i$. Since $j < m < n$, we have that $c_{n-1} = d_{n-1} = r$. 
     Then, $[d \harp m]_i = [e \harp m]_i$, so that $[f(d \harp m)]_i = [f(e \harp m)]_i$. 
     Hence, by (4), $[g(d)]_i =[g_{r+1}(d)]_i = [g_{r+1}(e)]_i= [g(e)]_i$.
     
     \smallskip
     
     $m \leq i < n$: Let $j = n-1$. Let $d,e \in \langle \omega \rangle^n$ be such that  $d_k = e_k$ whenever $ k < n-1$, intending to prove that $[g(d)]_i = [g(e)]_i$. If $d = e$, then the conclusion is trivial, so 
     suppose that $d_{n-1} = r < s = e_{n-1}$.  Let $c = d \harp m = e \harp m$. 
     Then,  $[f(c)]_i = [a_{c,i}]_i$. 
Therefore, by (3), we have that   $[g(d)]_i = [g_{r+1}(d)]_i = [a_{c,i}]_i = [g_{r+1}(e)]_i = [g(e)]_i$.

\smallskip

The claim $(*)$ is proved. By $n$ applications of Lemma~4.1, we get an infinite $Y \subseteq \omega$ such that for each $i < n$ there is  $S_i \subseteq n$ such that whenever $x,y \in \langle Y \rangle^n$, then $[g(x)]_i = [g(y)]_i$ iff 
$\{j<n : x_j \neq y_j\} \subseteq S_i$.  It follows from $(*)$ that each $S_i \neq \varnothing$. (In fact, from the proof of $(*)$, we get that $n \backslash m \subseteq S_i$ if $i < m$ and that $R_i \subseteq S_i$ if $m \leq i < n$.) Let $\vS = \langle S_0,S_1, \ldots, S_{n-1} \rangle$, and assume, without loss,  that $Y = \omega$. Then $g$ is an embedding of the $\vS$-halfcube into ${\mathcal A}$, which is a contradiction. \qed

\bigskip

It follows from Lemma~3.7 that if $\vS$ is an $n$-tuple of nonempty subsets of $m$ and $d = \tau(\vS)$,
then the $d$-halfcube is parbeddable into the $\vS$-halfcube. Therefore, the following corollary to Theorem~4.2 ensues.

\bigskip

{\sc Corollary 4.3}:  {\em Suppose that  ${\mathcal A}$ is an $n$-indexed hyperspace that does not parbed any $d$-halfcube, where $d \leq n$.              Then ${\mathcal A}$ has an acceptable coloring.}    

\bigskip

Restricting the previous corollary to $n$-grids, we get the following corollary that is a strengthening of de la Vega's Theorem~0.3.  

\bigskip

{\sc Corollary 4.4}: {\em Suppose that  ${\mathcal A}$ is an $n$-grid that does not  embed the $n$-halfcube. 
             Then ${\mathcal A}$ has an acceptable coloring.}    
             
             \bigskip
             
             {\it Proof}. Suppose that ${\mathcal A}$ is an $n$-grid. Then the only $d$-halfcube, where $d \leq n$, that it can parbed, is the $n$-halfcube. Any parbedding of the $n$-halfcube into ${\mathcal A}$ is a weak embedding. Finally, if ${\mathcal A}$ weakly embeds the $n$-halfube, then it embeds the $n$-halfcube. Thus, if ${\mathcal A}$ does not embed the $n$-halfcube, then it satisfies the hypothesis of Corollary~4.3. \qed 
             
             \bigskip

    Notice that Theorem~4.2 results when the hypothesis $d < \omega$ of the next theorem is replaced by $d = \infty$. 

             \bigskip


        {\sc Theorem 4.5}: {\em Suppose that  $1 \leq d < \omega$ and ${\mathcal A}$ is an $n$-indexed hyperspace  that  does not embed any $\vS$-halfcube, where      $\vS$ is an $n$-tuple of   subsets of $n$ and $\tau(\vS) < d$.   If $|A| < \aleph_{d-1}$, then ${\mathcal A}$ has an acceptable coloring.}    
    
    \bigskip

    {\it Proof}. This proof follows very closely the proof of Theorem~4.2. Theorem~2.4 gets used rather than Theorem~2.2.  There is one additional point 
    that needs to be checked. In the proof, we are assuming that $d < \omega$ and that $\vR$ 
    is an $I$-tuple of  nonempty subsets of $I$. It then must be shown that $\tau(\vR) >d-1$.     We then obtained the $n$-tuple  $\vS$ of nonempty 
    subsets of $n$ such that every finite $\vS$-cube is embeddable in ${\mathcal A}$. 
    This implies that $\tau(\vS) >d$. Thus, it remains to prove that 
    $\tau(\vR) \geq \tau(\vS) -1$.  But this is clear since if $T$ is a transversal of $\vR$ and $i \in I$, 
    then $T \cup \{i\}$ is a transversal of $\vS$.
    \qed
    
    \bigskip
    
        {\sc Corollary 4.6}: {\em Suppose that  $1 \leq d < n$ and  ${\mathcal A}$ is an $n$-indexed hyperspace  that  does not parbed the $(d-1)$-halfcube.    If $|A| < \aleph_{d-1}$,   then ${\mathcal A}$  has an acceptable coloring.}    
    
    \bigskip

    Suppose that ${\mathcal A}$ in Corollary~4.6 is an $n$-grid and $d = n-1$.   
    Since the $(d-1)$-halfcube is not  parbeddable into ${\mathcal A}$, then  Theorem~0.2 vacuously follows.

    \bigskip
    
   {\sc Definition 4.7}:  If ${\mathcal A}$ is an $n$-indexed hyperspace ${\mathcal A}$, then $\fcn({\mathcal A})$, the {\bf finite cube number} of ${\mathcal A}$, is the least $d$, where $1 \leq d \leq n$, such that for some $n$-tuple $\vS$ of subsets of $d$, ${\mathcal A}$ embeds every finite $\vS$-cube.  If there is no such $d$, then we let $\fcn({\mathcal A}) = \infty$. 
   
   \bigskip
    
    With this definition, we get the following corollary to Theorems~4.2 and~4.5.
    
    \bigskip
    
    {\sc Corollary 4.8}: {\em Suppose that ${\mathcal A}$ is an $n$-indexed hyperspace, $\fcn({\mathcal A})$ $ = d$ and $|A| < \aleph_{d-1}$. Then ${\mathcal A}$ has an acceptable coloring.} \qed
    
    \bigskip
    
    This corollary will be improved for semialgebraic indexed hyperspaces in the next section.
    
    
       \bigskip
       
        {\bf \S5.~Semialgebraic indexed hyperspaces.} Consider the ordered real field ${\widetilde \RR} = (\RR,+, \cdot, 0,1, \leq)$. We let ${\mathcal L}_{OF}$ be the language for $\widetilde \RR$. In this section, we will make tacit use of the famous theorems of Tarski that $\Th(\widetilde \RR)$, the first-order theory of ${\widetilde \RR}$, is decidable and admits the effective elimination of quantifiers.  If ${\widetilde R}$ is any ${\mathcal L}_{OF}$-structure and $X \subseteq R$, then ${\mathcal L}_{OF}(X)$ is ${\mathcal L}_{OF}$ augmented with (constants denoting) the elements of $X$. A subset $X \subseteq \RR^m$ is {\bf semialgebraic} if it is definable in ${\widetilde \RR}$ by a formula in which parameters are allowed. An $n$-indexed hyperspace ${\mathcal A} = (A;E_0,E_1, \ldots, E_{n-1})$ is {\bf semialgebraic} if, for some $m < \omega$, $A \subseteq \RR^m$ is semialgebraic as are each $E_i \subseteq \RR^{2m}$. If $\vS$ is an $n$-tuple of finite subsets of $m < \omega$, then the $\vS$-cube over $\RR$ is semialgebraic. Also, each $n$-indexed hyperspace 
        $(\RR^m;E(c_0),E(c_1), \ldots,$ $ E(c_{n-1}))$ from the prologue is semialgebraic. The purpose of this section is to generalize Theorem~0.4 from $n$-grids to $n$-indexed hyperspaces. 
        
        If $Y \subseteq  X_0 \times X_1 \times \cdots \times X_{m-1}$ and $f$ is a function on $Y$, 
        then $f$ is {\bf one-to-one in each coordinate} if whenever $x,y \in Y$, $i < m$ and 
        $x_j = y_j$ whenever $i \neq j < m$, then $f(x) = f(y) \Longleftrightarrow x = y$.  
        The following definition is adapted from \cite{grid}. 
        
        \bigskip
        
       {\sc Definition 5.1}:  Suppose that $\vS$ is an $n$-tuple  of subsets of $m < \omega$, ${\mathcal A}$ is an $n$-indexed hyperspace, and $X = X_0 \times X_1 \times \cdots \times X_{m-1}$.  A function $g : X \into A$ 
        is an  {\bf immersion} of the ${\vS}$-cube for $X$ into ${\mathcal A}$ if   the following hold:
        
        \begin{itemize}
        
        \item $g$ one-to-one in each coordinate; 
        
        \item if $x,y \in X$, $i<n$ and $g(x) \neq g(y)$, then \\ $[x]_i = [y]_i \Longleftrightarrow [g(x)]_i = [g(y)]_i$. 
        
              \end{itemize} 
        If there is an immersion of the $\vS$-cube for $X$ into ${\mathcal A}$, then we say that the $\vS$-cube for $X$ is {\bf immersible} into ${\mathcal A}$.  If $X = \RR^m$ and ${\mathcal A}$ is semialgebraic, then we say that the $\vS$-cube over $\RR$ is {\bf semialgebraically  immersible} into ${\mathcal A}$
        if there is a semialgebraic immersion $g : \RR^m \into A$.

        \bigskip
        
        If $g : X \into A$, where $X, \vS$ and ${\mathcal A}$ are as in Definition~5.1, then $g$ is an embedding of the $\vS$-cube for $X$ into ${\mathcal A}$ iff it is a one-to-one immersion.
        
        \bigskip
        
        {\sc Lemma 5.2}: {\em Let $\vS$ be an $n$-tuple of subsets of $m$ and ${\mathcal A}$  a semialgebraic $n$-indexed hyperspace. If the $\vS$-cube over $\RR$ is semialgebraically embeddable into ${\mathcal A}$, then there is a semialgebraic {\em analytic} embedding of the $\vS$-cube over $\RR$ into ${\mathcal A}$.}
        
        \bigskip
        
        {\it Proof}. Suppose that $f : \RR^m \into A$ is a semialgebraic embedding of the $\vS$-cube over $\RR$ into ${\mathcal A}$. By analytic cylindrical decomposition, there are disjoint analytic cylinders 
       $ B_0,B_1,\ldots,B_k \subseteq \RR^m$ whose union is $\RR^m$ and $f$ is analytic on each $B_i$.
        There is some $i \leq k$ such that $\dim(B_i) = m$. There are rationals $p_j < q_j$, for $j < m$, such that $B = (p_0,q_0) \times (p_1,q_1) \times \cdots \times (p_{m-1},q_{m-1}) \subseteq B_i$.
        Let $g_j : \RR \into (p_j,q_j)$ be an analytic, semialgebraic bijection, and let $g = (g_0,g_1, \ldots, g_{m-1})$. Then, $fg$ is a semialgebraic analytic embedding of the $\vS$-cube over $\RR$ into~${\mathcal A}$. \qed
        
        \bigskip

 We say that an $n$-tuple $\vS$ of subsets of $d$ is {\bf reduced} if $\tau(\vS) = d < \omega$.

        \bigskip
        
        {\sc Lemma 5.3}: {\em Suppose that ${\mathcal A}$ is a semialgebraic $n$-indexed hyperspace,  
        $\vS$ is a reduced $n$-tuple of subsets of $d$, and  the $\vS$-cube over $\RR$ is semialgebraically immersible into ${\mathcal A}$.  Then the $\vS$-cube over $\RR$ is embeddable into ${\mathcal A}$.}
        
        \bigskip
        
        {\it Proof}.  Let $f : \RR^d \into A$ be a semialgebraic immersion of the $\vS$-cube over $\RR$ into ${\mathcal A}$.   Let $\FF \subseteq \RR$ be a countable, real-closed subfield such that ${\mathcal A}$ is ${\mathcal \FF}$-semialgebraic and $f$ is $\FF$-definable.  Let $T$ be a transcendence basis for $T$ over $\FF$ such that whenever $a < b \in \RR$, then $|T \cap (a,b)| = 2^{\aleph_0}$.  For $i <d$, let $T_i = (i,i+1) \cap T$. Each $|T_i| = 2^{\aleph_0}$, so we have that   the $\vS$-cube for $T_0 \times T_1 \times \cdots \times T_{d-1}$ is isomorphic to the $\vS$-cube over $\RR$. We prove $(3)$ by proving  that 
$f \harp (T_0 \times T_1 \times \cdots \times T_{d-1})$ is an embedding of the $\vS$-cube
for $T_0 \times T_1 \times \cdots \times T_{d-1}$ into ${\mathcal A}$. Clearly, it suffices to prove that 
$f$ is one-to-one on $T_0 \times T_1 \times \cdots \times T_{d-1}$. 

For a contradiction, suppose that $s,t \in 
 T_0 \times T_1 \times \cdots \times T_{d-1}$, $s \neq t$ and   $f(s) = f(t)$.  Suppose that $i < d$ is such that $s_i \neq t$. For each $x \in \RR$, let $r(x) \in \RR^d$ be such that $r(x)_i = x$ and $r(x)$ agrees with $t$ on all other coordinates. Since $f$ is one-to-one on each coordinate, $f(s) = f(r(x))$ iff $x = t_i$. 
 But this gives an $\FF \cup (T \backslash \{t_i\})$-definition of $t_i$, contradicting that $T$ is algebraically independent over $\FF$.   \qed

 \bigskip
 
 {\sc Lemma 5.4}: {\em Suppose that ${\mathcal A}$ is a semialgebraic $n$-indexed hyperspace. 
 Then there is a finite partition $\{A_0,A_1, \ldots,A_m\}$ of $A$ such that for every $j \leq m$,
 there are $\ell < \omega$, an analytic semialgebraic bijection $f : A_j \into \RR^\ell$ and 
 analytic semialgebraic functions $e_0,e_1, \ldots,e_{n-1} : \RR^k \into \RR^k$ such that 
 for every $a,b \in A_j$ and $i < n$, $[a]_i = [b]_i$ iff $e_i(f(a)) = e_i(f(b))$.}
 
 \bigskip
 
 {\it Proof.} We give a sketch of the proof.  Let ${\mathcal A} = (A;E_0,E_1, \ldots, E_{n-1})$ be an $n$-indexed hyperspace where $A \subseteq R^k$. We are trying to get a finite partition ${\mathcal P}$ of $A$ as described in the lemma. 
 Let $g_0,g_1, \ldots, g_1 : A \into \RR^k$ be semialgebraic functions such that whenever $a,b \in A$ and $i < n$, then $[a]_i = [b]_i$ iff $g_i(a) = g_i(b)$. 
 Using analytic cylindrical cell decomposition, we get a semialgebraic partition $A = C_0 \cup C_1 \cup \cdots \cup C_t$ 
 such that each $C_j$ is an analytic cell and each $g_i$ is analytic on $C_j$. If $\dim(C_j) = \dim(A)$, then  put $C_j$ into ${\mathcal P}$.   Repeat process for each $C_j$ such that $\dim(C_j) < \dim(A)$.
 Continue putting cells into ${\mathcal P}$ until ${\mathcal P}$ is a partition of $A$ into cells 
 $A_0,A_1, \ldots, A_m$. For each $A_j$, there are $\ell \leq k$ and    an analytic semialgebraic bijection $f : A_j \into \RR^\ell$. \qed
 
 \bigskip
 
 The following theorem, which we refer to as the Polarized Canonical Erd\H{o}s-Rado Theorem (PCERT),
 will be needed. For more on this theorem, see, for example, \cite[Coro.~1.4]{sch16}). 
 If $X = X_0 \times X_1 \times \cdots \times X_{r-1}$ and $J \subseteq r$, then $\sim_J$ is the equivalence relation on $X$ induced by $J$; that is, if $x,y \in X$, then $x \sim_J j$ iff $x_i = y_i$ for all $i \in r \backslash J$. 
 
 \bigskip

{\sc Theorem 5.5}: (PCERT) 
{\em If $\lambda$ is a cardinal and $r < \omega$, then there is a cardinal $\kappa$ such that whenever  $\approx$ is an equivalence relation on $\kappa^r$, 
 then there are $J \subseteq r$ and $X_0,X_1, \ldots,X_{r-1} \subseteq \kappa$ such that $|X_0| = |X_1| = \cdots = |X_{r-1}| = \lambda$ and $\approx$ agrees with $\sim_J$ on  $X_0 \times X_1 \times \cdots \times X_{r-1}$.}

 \bigskip

{\sc Lemma 5.6}: {\em Suppose that ${\mathcal A}$ is a semialgebraic $n$-indexed hyperspace, 
  $\vS$ is a reduced $n$-tuple of subsets of $d< \omega$, and   every finite $\vS$-cube is embeddable into ${\mathcal A}$.  Then  the  
$\vS$-cube over $\RR$ is semialgebraically immersible into ${\mathcal A}$.}

\bigskip

{\it Proof}.    Let $\mathcal A$, $n$ and $\vS$   be as given.  Let $\{A_0,A_1, \ldots, A_m\}$ be a partition of $A$ as in Lemma~5.4. Since every finite $\vS$-cube is embeddable into ${\mathcal A}$, then (by  Finite Polarized Ramsey's Theorem) there is  $j \leq m$ such that every finite $\vS$-cube is embeddable into ${\mathcal A}|A_j$. Thus, we might as well assume that $A_j = A$. Then, using the function $f$ in Lemma~5.4, assume that $A = \RR^k$. Thus, we have   ${\mathcal A} = ( \RR^k; E_0,E_1, \ldots, E_{n_1})$, where $1 \leq k < \omega$, and analytic semialgebraic functions 
$e_0,e_1, \ldots, e_{n-1} : \RR^k \into \RR^k$ such that for each $i < n$ and $a,b \in \RR^k$, 
$[a]_i = [b]_i$ iff $e_i(a) = e_i(b)$.

\smallskip

Let $\widetilde R \succ \widetilde \RR$ be a sufficiently saturated elementary extension.  If $j < \omega$ and $D \subseteq \RR^j$, let $D^R$ be the subset of $R^j$ defined in ${\widetilde R}$ by a same formula that defines $D$ in ${\widetilde \RR}$. Let ${\mathcal A}^R = (R^d;E_0^R, E_1^R, \ldots, E_{n-1}^R)$, which is an $n$-indexed hyperspace. 
If $D \subseteq R^m$ and $X \subseteq R$, then we say that $D$ is $X$-definable if it is definable in ${\widetilde R}$ using only parameters from $X$.

Let $\FF \subseteq R$ be a countable real-closed subfield such that ${\mathcal A}^R$ and all the  $e_i^R$'s  are $\FF$-definable. 
Let $T \subseteq R$ be a transcendence basis for $\widetilde R$ over $\FF$ such that whenever $a,b \in \FF$ and $a < b$, then $|(a,b) \cap T| = |R|$. This choices of $\FF$ and $T$ are not definitive in that at various times in this proof we may replace $\FF$ by a larger real-closed field that is generated over $\FF$ by some finite subset $T_0 \subseteq T$. When we do that, it should be understood that we then replace $T$ by $T \backslash T_0$. 

If  $D \subseteq R^m$ is $R$-definable, then define $\supp(D)$, the {\it support} of $D$,  to be the smallest subset $S \subseteq T$ such that  $D$ is  $(S \cup \FF)$-definable. For each $R$-definable $D \subseteq R^m$,  $\supp(D)$ is a unique, finite subset of $T$. If $a \in R$ or $a \in R^k$, then $\supp(a) = \supp(\{a\})$.  If $a \in A$ and $i < n$, then $\supp([a]_i) \subseteq \supp(a)$.

Suppose that $1 \leq j < \omega$, $a \in R^j$ and  
$\supp(a) = \{t_0,t_1, \ldots, t_{m-1}\}_<$. (This notation implies that $t_0 < t_1 < \cdots < t_{m-1}$.)  A {\em determining function} for $a$ is an $\FF$-definable, ${\widetilde R}$-analytic function $f : \dom(f) \into R^j$  such that: 

\begin{itemize}

\item[(1)] $\dom(f)$ is an open subset of  $\langle R \rangle ^m$. (Recall that $\langle R \rangle ^m = \{x \in R^m : x_0 < x_1 < \cdots < x_{m-1}\}$.) 

\item[(2)] $\dom(f)$ is orthogonally convex (i.e., if $\ell \subseteq R^m$ is a line parallel to a coordinate axis, then $\ell \cap \dom(f)$ is convex). 

\item[(3)] $f$ is one-to-one in each coordinate. 

\item[(4)]  
$\langle t_0,t_1, \ldots, t_{m-1} \rangle \in \dom(f)$  and $f(t_0,t_1, \ldots, t_{m-1}) = a$.

\end{itemize}

\smallskip

{\sf Claim 1}: {\em Every $a \in R^j$ has a determining function.}

\smallskip

We sketch a proof since this is probably well known and, if not,  then the proof of a very similar statement (within the proof of \cite[Theorem~3.1]{sch16}) can be consulted. First, assume that $j=1$ so that $a \in R$. Let $\supp(a) = \{t_0,t_1, \ldots, t_{m-1}\}_<$.  Let $p(x,y_0,y_1, \ldots, y_{m-1}) \in R[x, \overline y]$ be such that 
$p(x,\overline t)$ is an irreducible polynomial and $p(a, \overline t) = 0$. Let $i < \omega$ be such that 
$a$ is the $i$-th root (in increasing order) of this polynomial. Then there is an $\FF$-definable function $g : D \into R$ such that 
$D \subseteq \langle R \rangle^m$, $\overline t \in D$ and $g(d)$ is the $i$-th root of $p(x, \overline d)$.
Using cylindrical cell decomposition for ${\widetilde R}$, we can get an $\FF$-definable, orthogonally convex  cell $C \subseteq D$ 
such that $t \in C$, $f = g \harp D$ is ${\widetilde R}$-analytic and 
$$\frac{\partial f}{\partial x_\ell}(d) \neq 0$$ whenever 
$d \in C$ and $\ell < m$. This $f \harp C$ is  a determining function for $a$..

Next, suppose that $j > 1$ and that $a \in R^j$. For $i < j$, let $f_i : C_i \into R$ be a determining function for $a_i$. These $f_i$'s can easily be merged into a function $f : C \into \RR^j$ that is a determining 
function for $a$. This completes the (sketch of) the proof of Claim~1. 

\smallskip

Let $g: X^d \into A^R$ be an embedding of the $\vS$-cube over $X$ into ${\mathcal A}^R$, were $X$ is sufficiently large. 
(It more than suffices to have $|X| \geq \beth_\omega$.) We will say that $Y \subseteq X^d$ is sufficiently large to mean that there are sufficiently large $X_0,X_1, \ldots, X_{d-1}$ such that 
$Y \supseteq X_0 \times X_1 \times \cdots \times X_{d-1}$
We can use PCERT  to get a sufficiently large $Y_0 \subseteq X^d$ such that:

\begin{itemize}

\item[(5)] There is a single $f$ that is a determining function for $g(x)$ whenever $x \in Y_0$. 

\end{itemize}

Let $m$ be such that $\dom(f) \subseteq R^m$. Notice that $m \geq 1$ since $|Y_0| \geq 2$. For each $x \in Y_0$, let 
$h(x) = \langle t_0,t_1, \ldots,t_{m-1} \rangle \in \dom(f)$, where $\supp(g(x)) = \{t_0,t_1, \ldots, t_{m-1}\}_<$. Thus, $h(x)_j$ is the $j$-th element in $\supp(g(x))$. Using PCERT again, we get a sufficiently large $Y_1 \subseteq Y_0$ such that:

\begin{itemize}

\item[(6)] Whenever $i \leq j < m$ and $x,y \in Y_1$, then $h(x)_i \leq h(y)_j$.

\end{itemize}

Thus, whenever $i < m$, then either for every $x,y \in Y_1$, then $h(x)_i = h(y)_i$ or else 
for every distinct $x,y \in Y_1$, then $h(x)_i = h(y)_i$. In the latter case, replace $\FF$ by the real-closed subfield of $R$ generated by $\FF$ and 
the common value $h(x)_i$. Thus, we can assume that $Y_1$  satisfies the following strengthening of $(6)$:

\begin{itemize}

\item[(6a)] Whenever $i < j < m$ and $x,y \in Y_1$, then $h(x)_i < h(y)_j$.

\item[(6b)] Whenever $i < m$ and $x,y \in Y_1$ are distinct, then $h(x)_i \neq h(y)_i$.

\end{itemize}

We next make a modification of $f$ and $\FF$. Because of (2),(3),(6b) and the saturation of $\widetilde R$, we can $r_0 < q_0 < r_1 < q_1 < \cdots < r_{m-1} < q_{m-1}$ in  $T$ such that:

\begin{itemize}

\item[(7)] Whenever $i < m$ and $x \in Y_1$, then $r_i < h(x)_i < q_i$. 

\item[(8)] $B = (r_0,q_0) \times (r_1,q_1) \times \cdots \times (r_{m-1}, q_{m-1}) \subseteq \dom(f)$.

\end{itemize}
We replace $\FF$ by its extension generated by $r_0,q_0,r_1,q_1, \ldots, r_{m-1},q_{m-1}$ and then  replace $f$ with $f \harp B$ so that we have

\begin{itemize}

\item[(9)] $\dom(f) = B$. 

\end{itemize}

Using PCERT again, we get a sufficiently large $Y_2 \subseteq Y_1$ such that:

\begin{itemize}

\item[(10)] For every $M \subseteq m$, there is $D_M \subseteq d$ such that whenever $x,y \in Y_2$,
then $x \sim_{D_M} y$ iff $h(x) \sim_M h(y)$.

\end{itemize}

\smallskip

{\sf Claim 2}: {\em If $M,N \subseteq m$, then $M \subseteq N$ iff $D_M \subseteq D_N$.}

\smallskip

We prove the claim. Consider $x,y \in Y_2$. Suppose $M \subseteq N$. Then  
$x \sim_{D_M} y$. Then, $x \sim_{D_M} y \Longrightarrow h(x) \sim_M h(y) \Longrightarrow h(x) \sim_N h(y) \Longrightarrow x \sim_{D_N} y$. This proves $M \subseteq N \Longrightarrow D_M \subseteq D_N$. For the converse, suppose that $D_M \subseteq D_N$. Then 
$h(x) \sim_{D_M} h(y) \Longrightarrow x \sim_M y \Longrightarrow x \sim_N y \Longrightarrow h(x) \sim_{D_N} h(y)$.

\smallskip

{\sf Claim 3}: {\em For each $i < n$, there is $M_i \subseteq m$ such that $D_{M_i} = S_i$.}

\smallskip

Fix $i < n$. Let 
$$
M_i = \{j < m : \exists x,y \in Y_2 \big([x]_i = [y]_i \wedge h(x)_j  \neq h(y)_j \big)\}.
$$

We first prove:

\smallskip

\noindent $(*)$ \hspace{30pt} $\forall j \in M_i \ \forall s,t \in B \cap T^d \big(s \sim_{\{j\}} t \into [f(s)]_i = [f(t)]_i\big)$. 

\smallskip

\noindent Let  $j \in M_i$. Let $x,y \in Y_2$ witness that $j \in M_i$. Let $s' = h(x)$ and $t' = h(y)$. Thus, $s'_j \neq t'_j$ and $[f(s')]_i = [f(t')]_i$. Then, $e^R_if(s') = e^R_if(t')$. 
Since $e^R_if$ is $R$-analytic and $\FF$-definable, it then follows that for every $s, t \in B \cap T^d$,
if $s \sim_{\{j\}} t$, then $e^R_if(s) = e^R_if(t)$, so that $[f(s)]_i = [f(t)]_i$. This proves $(*)$. 

We now prove that $D_{M_i} = S_i$. 

\smallskip

$D_{M_i} \subseteq S_i$: Suppose that $x \sim_{D_{M_i}} y$ (intending to show that $x \sim_{S_i} y$). 
Then, $h(x) \sim_{M_i} h(y)$. Let $t_0,t_1, \ldots, t_r \in B \cap T^d$ such that 
$t_0 = h(x)$, $t_r = h(y)$ and for all $\ell < r$ there is $j \in M_i$ such that $t_\ell \sim_{\{j\}} t_{\ell + 1}$. 
It follows from $(*)$ that $[g(x)]_i = [fh(x)]_i = [fh(y)]_i = [g(y)]_i$ so that $x \sim_{S_i} y$.

\smallskip

$S_i \subseteq D_{M_i}$: Suppose that $x \sim_{S_i} y$ (intending to show that $x \sim_{D_{M_i}} y$).
Then, $[x]_i = [y]_i$ so that  $h(x) \sim_{M_i} h(y)$ by the definition of $M_i$. Therefore, 
$x \sim_{D_{M_i}} y$. 

\smallskip

This completes the proof of Claim~3.   

\smallskip

We make two more modifications of $f$ and $\FF$. For the first one, suppose that there are $t\in T$ and  $j < m$ such that $h(x)_j = t$ whenever $x \in Y_2$. Replace $\FF$ by its extension generated by $t$ and then replace $f$ by the function $(m-1)$-ary function by fixing the $j$-th coordinate at $t$. We then have:

\begin{itemize}

\item [(11)] If $M \subseteq m$ and $D_M = \varnothing$, then $M = \varnothing$. 

\end{itemize}

Letting $M_i$ be as in Claim~3, it follows from Claim~2 and $(11)$, that $\tau(\langle M_0,M_2,\ldots, M_{d-1} \rangle) = d$. Let $I \subseteq m$ be a transversal for  $\langle M_0,M_1,\ldots, M_{d-1} \rangle$ 
such that $|I| = d$, where $I = \{i_0,i_1, \ldots, i_{d-1}\}_<$.  Let $t \in B \cap T^m$. We modify $f$ and $\FF$ by replacing $\FF$ with its extension generated by $\{t_j : j \in m \backslash I\}$. Let 
$B' = \{a \in B : a_j = t_j\}$ and then replacing $f$ by the function

 \bigskip
           
        {\sc Theorem  5.7}: {\em Suppose that ${\mathcal A}$ is a semialgebraic $n$-indexed hyperspace 
        and $d \leq n$.  The following are equivalent$:$
        
        \begin{itemize}
        
      \item[(1)] There is an $n$-tuple $\vS$ of subsets of $d$ such that  the $\vS$-cube over $\RR$ is semialgebraically immersible into ${\mathcal A}$.
       
        \item[(2)] There is an $n$-tuple $\vS$ of subsets of $d$ such that  every finite $\vS$-cube is embedable into ${\mathcal A}$.
        
        \item[(3)] There is an  $n$-tuple $\vS$ of subsets of $d$ such that  the $\vS$-cube over $\RR$ is embeddable into ${\mathcal A}$.
        
        \end{itemize}}
        
        \bigskip
        
        {\it Proof}. $(3) \Longrightarrow (2)$ is trivial. Lemma~5.6 implies $(2) \Longrightarrow (1)$ and Lemma~5.3 implies $(1) \Longrightarrow (3)$. \qed
        
        \bigskip
        
        If ${\mathcal A}$ is a semialgebraic $n$-indexed hyperspace, then $\fcn({\mathcal A})$ (see Definition~4.7) is the least $d$ $(1 \leq d \leq n)$ such that every (or any) one of $(1)$ -- $(3)$ holds. If there is no such $d$, then 
        $\fcn({\mathcal A}) = \infty$.

           \bigskip
    
    {\sc Corollary 5.8}: {\em Suppose that ${\mathcal A}$ is a semialgebraic $n$-indexed hyperspace and  $\fcn({\mathcal A}) = d$. Then ${\mathcal A}$ has an acceptable coloring iff $2^{\aleph_0} < \aleph_{d-1}$.} \qed
 
\bigskip

{\sc Corollary 5.9}: {\em The set of  ${\mathcal L}_{OF}$-formulas  that, for some $n < \omega$, define in ${\widetilde {\RR}}$ a semialgebraic $n$-indexed hyperspace having an acceptable coloring is computable.}

\bigskip

{\it Proof}. Let $\Gamma$ by the set of ${\mathcal L}_{OF}$-formulas defined in the corollary. 
Using $(3) \Longrightarrow (1)$ of the applicable  one of Corollary~5.3 or 5.6, we get that $\Gamma$ is c.e., and using  $(3) \Longleftrightarrow (2)$ we get that $\Gamma$ is co-c.e. \qed  

\bigskip

In the previous corollary, the formulas are ${\mathcal L}_{OF}$-formulas, so they are not allowed to have any parameters. There 
is a way to modify this corollary for ${\mathcal L}_{OF}(\RR)$-formulas. A typical ${\mathcal L}_{OF}(\RR)$-formula has the form $\varphi(x,c)$, where  $\varphi(x,y)$ is an $(m+n)$-ary ${\mathcal L}_{OF}$-formula and $c \in \RR^n$. We say that a  set $\Gamma$  of ${\mathcal L}_{OF}(\RR)$-formulas $\varphi(x,c)$ 
 is {\bf decidable} if there is a computable set $\Delta$ of ${\mathcal L}_{OF}$-formulas such that for every ${\mathcal L}_{OF}(\RR)$-formula $\varphi(x,c)$, the following are equivalent:

\begin{itemize}

\item [(1)] $\varphi(x,c) \in \Gamma$;

\item[(2)] there is a formula $\theta(y) \in \Delta$ such that 
${\widetilde \RR} \models \theta(c)$ and \\ $\forall y[\theta(y) \into \varphi(x,y)]$ is in $\Delta$;

\item[(3)] there is a formula $\theta(y) \in \Delta$ such that 
${\widetilde \RR} \models \theta(c)$ and \\ $\forall y[\theta(y) \into \neg\varphi(x,y)]$ is in $\Delta$.

\end{itemize}
A set of ${\mathcal L}_{OF}$-formulas is computable iff it is decidable (as a set of ${\mathcal L}_{OF}(\RR)$-formulas). 

\bigskip

{\sc Corollary 5.10}: {\em The set of  ${\mathcal L}_{OF}(\RR)$-formulas  that define in $\widetilde{\RR}$ a semialgebraic indexed hyperspace having an acceptable coloring is decidable.}
\qed


\bibliographystyle{plain}

\end{document}